\documentclass{amsart}
\usepackage{amssymb,latexsym,epsfig,color}
 
\newtheorem{theorem}{Theorem}[section]
\newtheorem{proposition}[theorem]{Proposition}
\newtheorem{corollary}[theorem]{Corollary}

\newtheorem{lemma}[theorem]{Lemma}

\newtheorem{problem}[theorem]{Problem}

\theoremstyle{definition}
\newtheorem{definition}[theorem]{Definition}
\newtheorem{remark}[theorem]{Remark}
\newtheorem{example}[theorem]{Example}

\def\endproof{\hfill$\square$\medskip}

\def\N{\ensuremath{\mathbb{N}}}
\def\Z{\ensuremath{\mathbb{Z}}}
\def\Q{\ensuremath{\mathbb{Q}}}
\def\R{\ensuremath{\mathbb{R}}}

\newcommand{\B}{{\mathcal B}}

\def\F{\ensuremath{\mathcal{F}}}

\def\A{\ensuremath{\mathcal{A}}}
\def\C{\ensuremath{\mathcal{C}}}

\def\a{\ensuremath{{\bf{a}}}}
\def\b{\ensuremath{{\bf{b}}}}
\def\cc{\ensuremath{{\bf{c}}}}
\def\d{\ensuremath{{\bf{d}}}}
\def\e{\ensuremath{{\bf{e}}}}
\def\k{\ensuremath{{\bf{k}}}}
\def\m{\ensuremath{{\bf{m}}}}
\def\p{\ensuremath{{\bf{p}}}}
\def\q{\ensuremath{{\bf{q}}}}

\def\u{\ensuremath{{\bf{u}}}}
\def\v{\ensuremath{{\bf{v}}}}
\def\x{\ensuremath{{\bf{x}}}}
\def\y{\ensuremath{{\bf{y}}}}
\def\z{\ensuremath{{\bf{z}}}}

\def\ica{\ensuremath{{I_{{\mathcal C}_{\mathcal A}}}}}
\def\ia{\ensuremath{{{I_{\mathcal A}}}}}

\def\inomega{\ensuremath{\textup{in}_{\omega}}}
\def\inu{\ensuremath{\textup{in}_{\u}}}
\def\inv{\ensuremath{\textup{in}_{\v}}}
\def\tin{\textup{in}}
\def\la{\ensuremath{{{\mathcal{L}_{\mathcal A}}}}}
\def\na{\ensuremath{{{\mathbb{N} A}}}}
\def\xa{\ensuremath{{{X_{\mathcal A}}}}}
\def\hia{\ensuremath{{{\mathbf{H}_{\ia}}}}}
\def\hica{\ensuremath{{{\mathbf{H}_{\ica}}}}}

\def\max{\ensuremath{\textup{max}}}
\def\min{\ensuremath{\textup{min}}}
\def\supp{\ensuremath{\textup{supp}}}
\def\cone{\ensuremath{\textup{cone}}}
\def\conv{\ensuremath{\textup{conv}}}
\def\deg{\ensuremath{\textup{deg}}}
\def\ker{\ensuremath{\textup{ker}}}
\def\ass{\ensuremath{\textup{Ass}}}
\def\top{\ensuremath{\textup{top}}}

\newcommand{\bigo}{{\mathcal O}}

\title [Circuit Ideals]
{The Circuit Ideal of a Vector Configuration}

\author{Tristram Bogart}
\address{Department of Mathematics, University of
  Washington, Seattle, WA 98195-4350}
\email{bogart@math.washington.edu}
\author{Anders N. Jensen}
\address{Institut for Matematiske Fag, Aarhus Universitet, DK-8000, 
{\AA}rhus, Denmark}
\email{ajensen@imf.au.dk}
\author{Rekha R. Thomas}
\address{Department of Mathematics, University of
  Washington, Seattle, WA 98195-4350}
\email{thomas@math.washington.edu}
\thanks{Tristram Bogart was partially supported by NSF grant
  DMS-0354131, Anders N.  Jensen by the Danish Research Training
  Council (Forskeruddannelsesr\aa det, FUR) and the Swiss National
  Science Foundation Project 200021-105202, and Rekha R. Thomas by NSF
  grant DMS-0401047.}  
\date{\today}

\begin{document}

\begin{abstract} 
  The circuit ideal, $\ica$, of a configuration $\A = \{\a_1, \ldots,
  \a_n\} \subset \Z^d$ is the ideal generated by the binomials
  ${\x}^{\cc^+} - {\x}^{\cc^-} \in \k[x_1, \ldots, x_n]$ as $\cc =
  \cc^+ - \cc^- \in \Z^n$ varies over the circuits of $\A$. This ideal
  is contained in the toric ideal, $\ia$, of $\A$ which has numerous
  applications and is nontrivial to compute.  Since circuits can be
  computed using linear algebra and the two ideals often coincide, it
  is worthwhile to understand when equality occurs.
  
  In this paper we study $\ica$ in relation to $\ia$ from various
  algebraic and combinatorial perspectives. We prove that the
  obstruction to equality of the ideals is the existence of certain
  polytopes. This result is based on a complete characterization of
  the standard pairs/associated primes of a monomial initial ideal of
  $\ica$ and their differences from those for the corresponding toric
  initial ideal.  Eisenbud and Sturmfels proved that $\ia$ is the
  unique minimal prime of $\ica$ and that the embedded primes of
  $\ica$ are indexed by certain faces of the cone spanned by $\A$. We
  provide a necessary condition for a particular face to index an
  embedded prime and a partial converse. Finally, we compare various
  polyhedral fans associated to $\ia$ and $\ica$. The Gr\"obner fan of
  $\ica$ is shown to refine that of $\ia$ when the codimension of the
  ideals is at most two.

\end{abstract}
\maketitle

\section{Introduction}
\label{sec:introduction}
Throughout this paper, we fix an ordered vector configuration  
$\A = \{\a_1, \ldots, \a_n \} \subset \Z^d$.  Assume that the $d \times n$
integer matrix $A = [\a_1 \, \ldots \,\a_n]$ whose columns are the
elements of $\A$ has rank $d$. Let $\la$ be the $(n-d)$-dimensional
saturated lattice $\{ \u \in \mathbb Z^n \, : \, A\u = {\bf 0} \}$.  
We assume that $\la \cap \N^n = \{ {\bf 0} \}$.

The {\em support} of a vector $\u \in \mathbb Z^n$ is defined to be
$\supp(\u) := \{ i \,:\, u_i \neq 0 \}$ and $\u$ is {\em primitive} if
the greatest common divisor of its components is one.

\begin{definition}
  A vector $\cc \in \la$ is a {\bf circuit} of $\A$
  if (1) $\cc$ is a non-zero primitive vector and (2) there does not
  exist $\d \in \la$ with $\supp(\d) \subsetneq \supp(\cc)$.
\end{definition}

Let $\C_{\A}$ denote the set of all circuits of $\A$. Write $\cc = \cc^+
- \cc^-$ where $c^+_j = c_j$ if $c_j > 0$ and $0$ otherwise, and
$c^-_j = -c_j$ if $c_j < 0$ and $0$ otherwise. Identify $\cc \in
\C_{\A}$ with the binomial $\x^{\cc^{+}} - \x^{\cc^{-}} \in \k[x_1,
\ldots, x_n] =: \k[\x]$ where $\k$ is an algebraically closed field and 
$\x^\u := x_1^{u_1} x_2^{u_2} \cdots x_n^{u_n}$.  We refer to both
$\cc$ and $\x^{\cc^{+}} - \x^{\cc^{-}}$ as a circuit of $\A$ and 
denote both lists by $\C_{\A}$.

\begin{definition}
  The {\bf circuit ideal} of $\A$ is the binomial ideal $ \ica :=
  \langle \C_{\A} \rangle \subseteq \k[\x].$
\end{definition}

The circuit ideal $ \ica $ is a subideal of the binomial prime 
\emph{toric ideal} of $\A$
$$
\ia := \langle \x^{\u^{+}} - \x^{\u^{-}} \, : \, \u \in \la
\rangle.$$
Toric ideals are the defining ideals of \emph{toric varieties}~\cite{Ful} 
and have numerous applications in combinatorics, optimization, algebra and
algebraic geometry~\cite{GBCP}. These connections make the computability 
of $\ia$ an important practical concern. 

\begin{proposition} \cite{GBCP} \label{prop:saturation}
  Given a finite subset $\B$ of $\la$, define the ideal $J_{\B} :=
  \langle \x^{\b^{+}} - \x^{\b^{-}} \, : \, \b \in \B \rangle \subset
  \k[\x]$. A set $\B$ spans $\la$ if and only if $(J_{\B} \, : \,
  (x_1x_2 \cdots x_n)^{\infty}) = \ia$.
\end{proposition}

Proposition~\ref{prop:saturation} is the starting point of the best
algorithms to compute $\ia$ since a spanning set $\B$ of $\la$ can be
computed easily and each saturation in
$$(J_{\B} \, : \, (x_1x_2 \cdots x_n)^{\infty}) = ((((J_{\B} \, : \,
x_1^{\infty}) \, : \, x_2^{\infty}) \cdots ) \, : \, x_n^{\infty})$$
can be achieved by a Gr\"obner basis calculation (\cite{HS},
\cite[Chapter 12]{GBCP}.)  It can be checked that ${\C}_{\A}$ spans
$\la$ and hence $\ia = (\ica \,:\, (x_1x_2 \cdots x_n)^{\infty})$.  In
many examples, $\ica$ equals $\ia$, and since the circuits of a matrix
can be computed easily~\cite[page 190]{DES}, it is of interest to know 
how close the circuit ideal is to the toric ideal and in particular when 
they are equal. This raises the main problem addressed in this paper. See 
Remark~\ref{rem:motivations} for further motivations.

\begin{problem} \label{equality}
When does the circuit ideal $\ica$ equal the toric ideal $\ia$?
\end{problem}

In this paper, we investigate Problem~\ref{equality} from several
different angles. Let $\na$ denote the semigroup $\{A\u \,:\, \u \in
\N^n \} \subset \Z^d$. Both $\ia$ and $\ica$ are homogeneous under
multi-grading by $\na$ with $\k[\x] / \ia$ having Hilbert function
value one for all $\b \in \na$.  In Section~\ref{sec:optimization} we
recall conditions for the equality of $\ia$ and $\ica$ and then
exhibit various properties of circuit ideals that contrast those of
toric ideals. We interpret the multi-graded Hilbert function values of
$\k[\x] / \ica$.

From the point of view of Gr\"obner basis theory, it is natural to
investigate $\ia$ and $\ica$ by examining the difference between their
initial ideals with respect to a fixed weight vector $\omega$.  In
Section~\ref{sec:initialideals}, we give a complete characterization
of the associated primes of a monomial initial ideal of $\ica$
(Theorem~\ref{thm:stdpairs}) extending previously known
characterizations of the associated primes of a monomial initial ideal
of $\ia$ \cite{HT2}. The associated primes and the difference between
the two monomial initial ideals are described in terms of certain
polytopes that depend on $\A$ and $\omega$. Using this we answer
Problem~\ref{equality} by showing that the obstruction to equality of
the ideals is the existence of certain polytopes of the above type 
(Theorem~\ref{thm:obstruction}).

In \cite{EiSt}, Eisenbud and Sturmfels showed that $\ia$ is the unique
minimal prime of $\ica$. Thus a second natural measure of the
difference between the two ideals is an understanding of the embedded
primes of $\ica$.  Let $\cone(\A)$ denote the $d$-dimensional cone
spanned by $\A$.  Record a face $\sigma$ of $\cone(\A)$ as the set of
indices, $j$, of all $\a_j$ that lie on $\sigma$.  Eisenbud and
Sturmfels proved that the associated primes of $\ica$ are all of the
form $P_{\sigma} + \ia$ where $\sigma$ is some face of $\cone(\A)$ and
$P_{\sigma}:= \langle x_j \, : \, j \not \in \sigma \rangle$.  In
particular, $\ia = P_{[n]} + \ia$ is indexed by the full face $[n] :=
\{1,2,\ldots,n\}$ of $\cone(\A)$. However, not all faces of
$\cone(\A)$ need index an associated prime of $\ica$ and Eisenbud and 
Sturmfels raise the following problem.

\begin{problem} \cite[\S 7]{EiSt} \label{prob:assprimes}
``It remains an interesting combinatorial problem to characterize the 
embedded primary components of the circuit ideal $\ica$.  In 
particular, which faces of $\cone(\A)$ support an associated prime 
of $\ica$?  An answer to this question might be valuable for 
the applications of binomial ideals to integer programming and 
statistics.''
\end{problem}

In Section~\ref{sec:assprimes}, we give a necessary condition for a
prime $P_{\sigma} + \ia$ to be an embedded prime of $\ica$
(Theorem~\ref{thm:primestoCSPs}) using the results in
Section~\ref{sec:initialideals}.  We also provide a partial converse
to Theorem~\ref{thm:primestoCSPs}. As an application, we derive
connections between the smoothness of the toric variety defined by a 
face $\sigma$ of $\cone(\A)$ and $P_{\sigma} + \ia$ being an
embedded prime of $\ica$ when $\A$ is a graded vector configuration.  

We conclude the paper in Section~\ref{sec:fans} by examining various
polyhedral fans associated to $\ia$ and $\ica$. Given a homogeneous
ideal $I$ and a weight vector $\omega \in \R^n$, let $\inomega(I)$ be
the initial ideal of $I$ with respect to $\omega$,
$\sqrt{\inomega(I)}$ the radical of $\inomega(I)$, and
$\top(\inomega(I))$ the intersection of the top-dimensional primary
components of $\inomega(I)$.  These entities define three equivalence
relations on $\R^n$ as follows.
\begin{enumerate}
\item The \emph{initial ideal} equivalence relation: $\u \sim \v
\Leftrightarrow \inu(I)=\inv(I)$, 
\item the \emph{top} equivalence relation: $ \u \sim \v
\Leftrightarrow \top(\inu(I))=\top(\inv(I))$, and 
\item the \emph{radical} equivalence relation: $\u\sim \v
\Leftrightarrow \sqrt{\inu(I)}=\sqrt{\inv(I)}$.
\end{enumerate}
For any homogeneous ideal $I$, the initial ideal equivalence classes
form the cells of the {\em Gr\"obner fan} of
$I$~\cite{MR},~\cite[Chapter 2]{GBCP}.  For $\ia$ it is well known
that the other two equivalence classes also form polyhedral fans ---
the radical equivalence relation gives the {\em secondary fan} of
$\A$~\cite{BFS},~\cite[Chapter 8]{GBCP}, and the top equivalence 
relation gives the \emph{hypergeometric fan} of $\A$~\cite{SST}.
We prove that for $\ica$, the top dimensional equivalence classes of
the radical and top equivalence relations coincide with those for
$\ia$ (Theorem~\ref{thm:circuitrefinement} and
Proposition~\ref{prop:top}).  However, the Gr\"obner fans of $\ia$ and
$\ica$ do not coincide in general. Corollary~\ref{cor:codim2} proves
that when the codimension of the ideals is at most two, the Gr\"obner
fan of $\ica$ refines that of $\ia$.

\section{Properties of $\ia$ versus $\ica$}
\label{sec:optimization}

In this section we first collect conditions equivalent to the
equality of $\ia$ and $\ica$. Many of these stem from combinatorics
and optimization and most are well known~\cite{DES},~\cite{GBCP}. 
We then contrast $\ica$ with $\ia$ in light of these conditions.

Consider the semigroup homomorphism $\pi: \, \N^n \rightarrow \na \, ,
\; \u \mapsto A\u$.  Both $\ia$ and $\ica$ are homogeneous under the
\emph{$\A$-grading} of $\k[\x]$ by $\deg(x_i) = \a_i$ for $i=1,
\ldots, n$ since every binomial of the form $\x^\u - \x^\v$ in either
ideal is \emph{$\A$-homogeneous} with \emph{$\A$-degree} $\pi(\u) = A
\u = A \v = \pi(\v)$.  Let $\hica \,:\, \na \rightarrow \N$ be the
$\A$-graded Hilbert function of $\k[\x]/\ica$ given by $\b \mapsto
\dim_{\k}(\k[\x]/\ica)_{\b}$.  Let $\hia$ be the same for $\ia$.

Since $\la \cap \N^n = \{ {\bf 0} \}$, for each $\b \in \na$, the
polyhedron $P_{\b} := \{ \x \in \R^n_{\geq 0} \,:\, A\x = \b \}$ is
bounded~\cite{Sch} which implies that $\pi^{-1}(\b) := \{ \x \in \N^n
\,:\, A\x = \b \} = P_{\b} \cap \N^n$ is finite for all $\b \in \na$.
For a fixed $\b \in \na$, the vertex set $\pi^{-1}(\b)$ admits four
natural graphs as follows. First, choose binomial generating sets
$G(\ia)$ and $G(\ica)$ of $\ia$ and $\ica$ respectively, and let
$\F(\b)$ and $\F^{\C}(\b)$ be graphs on $\pi^{-1}(\b)$ such that $\u$
is adjacent to $\v$ in $\F(\b)$ (respectively in $\F^{\C}(\b)$) if
$\x^\u-\x^\v$ is a monomial multiple of a binomial in $G(\ia)$
(respectively in $G(\ica)$).  Next, fix a generic weight vector
$\omega \in \R^n$ in the sense that the initial ideals $\inomega(\ia)$
and $\inomega(\ica)$ are both monomial ideals.  Let $G_{\omega}(\ia)$
and $G_{\omega}(\ica)$ be the marked reduced Gr\"obner bases of $\ia$
and $\ica$ with respect to $\omega$.  Elements of these Gr\"obner
bases are $\A$-homogeneous binomials and the Gr\"obner basis being
marked means that the first term in each binomial $f$ is its initial
term $\inomega(f)$. Construct the \emph{directed} graphs
$\F_{\omega}(\b)$ and $\F^{\C}_{\omega}(\b)$ on $\pi^{-1}(\b)$ by
drawing an arrow from $\u$ to $\v$ in $\F_{\omega}(\b)$ (respectively
$F^{\C}_{\omega}(\b)$) if and only if $\x^{\u} - \x^{\v}$ is a
monomial multiple of some marked binomial in $G_{\omega}(\ia)$
(respectively $G_{\omega}(\ica)$). 

\begin{lemma} \cite[Theorem~1.1]{DES} 
\label{lem:samecomponent}
Vectors $\u, \v \in \pi^{-1}(\b)$ are in the same component of
$\F(\b)$ (respectively $\F^{\C}(\b)$) if and only if $\x^{\u} -
\x^{\v}$ lies in $\ia$ (respectively $\ica$).
\end{lemma}

Lemma~\ref{lem:samecomponent} shows that while the edges in $\F(\b)$
and $\F^{\C}(\b)$ depend on the choice of generating sets $G(\ia)$ and
$G(\ica)$, the components, and in particular the number of components,
do not depend on this choice. Further, both $\F(\b)$ and
$\F_\omega(\b)$ partition $\pi^{-1}(\b)$ identically into components.
The same holds for $\F^{\C}(\b)$ and $\F_\omega^{\C}(\b)$.  The
following theorem mostly collects results from~\cite{DES} 
and~\cite{GBCP}.
 
\begin{theorem} 
\label{thm:conditionsforequality}
The following statements are equivalent.
\begin{enumerate}
\item The ideals $\ia$ and $\ica$ are equal.
\item For every $\b \in \na$, the graph $\F^{\C}(\b)$ is connected.
\item For every $\b \in \na$, the digraph $\F^{\C}_\omega(\b)$ has a
  unique sink. (In this case, the unique sink $\u$ in
  $\F^{\C}_\omega(\b)$ is the optimal solution of the integer program
  $ \textup{minimize} \, \{ \omega \cdot \x \,:\, A \x = \b, \, \x \in
  \N^n \}$.)
\item For every $\b \in \na$ and generic weight vector $\omega \in
  \R^n$, $\inomega(\ica)$ has a unique standard monomial of
  $\A$-degree $\b$. (In this case, the standard monomial of
  $\A$-degree $\b$ is $\x^\u$ where $\u$ is the unique sink in
  $\F^{\C}_\omega(\b)$.)
\item For every $\b \in \na$, the Hilbert function value $\hica(\b)$
  is one.
\end{enumerate}
\end{theorem}

\begin{proof}
  Statements (2)--(5) are all true if $\ica$ is replaced by $\ia$,
  $\F^{\C}(\b)$ by $\F(\b)$ and $\F^{\C}_\omega(\b)$ by
  $\F_\omega(\b)$, see \cite[Chapters 4,5,10]{GBCP}.  Further, $\ia =
  \ica$ if and only if $G_\omega(\ia) = G_\omega(\ica)$ if and only if
  for each $\b$, $\F^{\C}_\omega(\b)$ equals $\F_\omega(\b)$ and hence
  if and only if $\F^{\C}(\b)$ and $\F(\b)$ have the same components.
  Hence (1) is equivalent to (2) and (3). Since $\ica \subseteq \ia$,
  the two ideals are equal if and only if (5). The equivalence of (3)
  and (4) follows from Lemma~\ref{lem:cantconnect} below.
\end{proof}

\begin{remark} \label{rem:motivations}
\begin{enumerate}
\item The connectivity of $\F(\b)$ was used in~\cite{DS98}, in the
  context of statistical sampling, to devise random walks on
  $\pi^{-1}(\b)$.  The equality of $\ia$ and $\ica$ allows
  $\pi^{-1}(\b)$ to be connected using $G(\ica)$, which is cheaper to
  compute than $G(\ia)$. In Section~\ref{sec:initialideals} we will
  see that for most $\b \in \na$, $\F^{\C}(\b)$ is in fact connected
  (Theorem~\ref{thm:hilbfunctionone}), and that the set of $\b \in
  \na$ for which $\F^{\C}(\b)$ is disconnected can be described
  precisely.  See also~\cite{DES}.

\item The equality of $\ia$ and $\ica$ will allow all integer programs of
the form
$$\textup{minimize} \, \{ \omega \cdot \x \,:\, A \x = \b, \x \in \N^n
\}$$
as $\b$ and $\omega$ vary to be solved by reduced Gr\"obner bases
of $\ica$. The significance of this is that the circuits of $\A$ are
precisely the primitive edge directions of the polyhedra $P_{\b}$ as
$\b$ varies in $\na$ and hence the directions taken by the simplex
algorithm in solving linear programs of the form
$$
\textup{minimize} \, \{ \omega \cdot \x \,:\, A \x = \b, \x \in
\R^n_{\geq 0} \}$$
as $\b$ and $\omega$ vary. Hence, even though
$\ica$ contains more than the circuits of $\A$, philosophically, the
equality of $\ia$ and $\ica$ allows \emph{integer} programs arising
from $A$ to be solved via the corresponding \emph{linear} programming
data.
\end{enumerate}
\end{remark} 

We now exhibit various properties of $\ica$ that contrast those of
$\ia$.

\begin{proposition}
\label{prop:manycomponents}
\begin{enumerate}
\item The graph $\F^{\C}(\b)$ may have arbitrarily many components,
  even if we restrict to the case of $A \in \Z^{1 \times 3}$.
\item The standard monomials of $\inomega(\ica)$ of  
$\A$-degree $\b$ are not necessarily the cheapest monomials of that 
degree with respect to $\omega$.
\end{enumerate}
\end{proposition}

\begin{proof}
(1) For any natural number $k \geq 2$, take $ A_k = \left(
    \begin{array}{rrr} k & 2k+1 & 3k+1 \end{array} \right).$ Since the
  three entries are pairwise relatively prime, the circuits are
  $x^{2k+1}-y^k$, $x^{3k+1}-z^k$, and $y^{3k+1}-z^{2k+1}$.  Their
  respective $\A$-degrees are $2k^2+k, 3k^2+k$, and $6k^2+5k+1$.  Thus
  the graph $\F^{\C}(b)$ has no edges when $b < 2k^2+k$.  In
  particular, this holds if we take $b = m(3k+1)$ for $m = \lfloor k/2
  \rfloor$.  This graph has at least $m+1$ vertices $\{(j,j,m-j): 0
  \leq j \leq m \}$, so it has at least $m+1 = \lfloor k/2 \rfloor +
  1$ components.
  
  (2) Consider $\A = \{3,4,5\}$.  The graded reverse lexicographic
  (grevlex) Gr\"obner basis of $\ica$ with $a \succ b \succ c$ is
  $$
  \{a^4-b^3, ab^3-c^3, b^5-c^4, b^2c^3-ac^4,
  a^3c^3-bc^4,a^2bc^4-c^6\}.$$
  The monomials of degree $17$ are $a^4c,
  a^3b^2, abc^2$ and $b^3c$ of which the last three are standard
  monomials of the above grevlex initial ideal of $\ica$. However, we
  see that the non-standard monomial $a^4c$ is cheaper than the
  standard monomial $a^3b^2$.
\end{proof}

We now prove that $\hica(\b)$ equals the number of components of
$\F^{\C}_\omega(\b)$, or equivalently, $\F^{\C}(\b)$.
Proposition~\ref{prop:manycomponents} (1) then shows that the values
of $\hica$ can be arbitrarily large even for $d$ and $n$ fixed. In
contrast, $\hia(\b) = 1$ for all $\b \in \na$.

\begin{lemma} 
\label{lem:cantconnect}
\begin{enumerate}
\item Let $\x^{\u}$ be a standard monomial of $\inomega(\ica)$ with $A\u =
\b$.  Then for all $\v \neq \u$ in the same component of
$\F^{\C}_{\omega}(\b)$ as $\u$, $\omega \cdot \u < \omega \cdot \v$.
In particular, if $\p \in \N^n$, $A\p = \b$, and $\omega \cdot \p <
\omega \cdot \u$, then $\x^{\u} - \x^{\p} \not \in \ica$. 
\item Each component of $\F^{\C}_{\omega}(\b)$ has a unique sink $\u$ and
  $\x^{\u}$ is the unique standard monomial of
  $\inomega(\ica)$ among all monomials $\x^{\v}$ such that
  $\v$ is in the same component as $\u$.
\end{enumerate}
\end{lemma}

\begin{proof} 
  (1) If $\v \,(\neq \u)$ lies in the same component of
  $\F^{\C}_{\omega}(\b)$ as $\u$, then by
  Lemma~\ref{lem:samecomponent}, $0 \neq f := \x^{\v} - \x^{\u} \in
  \ica$. Since $\x^{\u} \not \in \inomega(\ica)$, $\inomega(f)$ equals
  either $\x^\v$ or $f$.  By the genericity of $\omega$, we can assume
  $\inomega(f) \neq f$. Thus $\omega \cdot \u < \omega \cdot \v$.
  
  (2) Let $D$ be an arbitrary component of $\F^{\C}_{\omega}(\b)$,
  $\v$ be an arbitrary vertex in $D$, and $\x^{\u}$ be the normal form
  of $\x^{\v}$ with respect to $G_{\omega}(\ica)$.  Then $\x^{\u}$ is
  a standard monomial of $\inomega(\ica)$ and by
  Lemma~\ref{lem:samecomponent}, $\u$ is in $D$. If $\x^{\u'}$ is
  another standard monomial of $\inomega(\ica)$ with $\u'$ in $D$,
  then by Lemma~\ref{lem:samecomponent}, $f := \x^\u - \x^{\u'} \in
  \ica$ with $\inomega(f)$ equal to either $\x^\u$ or $\x^{\u'}$, a
  contradiction. This implies that $\x^\u$ is the unique normal form
  of all $\x^\v$, $\v \in D$ and hence it is the unique sink in $D$.
\end{proof}

\begin{proposition} \label{prop:countcomponents} 
  The Hilbert function value $\hica(\b)$ equals the number of
  components of $\F^{\C}_{\omega}(\b)$.
\end{proposition} 

\begin{proof}
  By Lemma~\ref{lem:cantconnect} (2), each component of
  $\F^{\C}_{\omega}(\b)$ contributes precisely one standard monomial
  of $\inomega(\ica)$. The number of standard monomials of
  $\inomega(\ica)$ of degree $\b$ equals
  $\dim_{\k}(\k[\x]/\inomega(\ica))_{\b} =
  \dim_{\k}(\k[\x]/\ica)_{\b} = \hica(\b)$.
\end{proof}

\begin{example} \label{hf picture}
When $\ia \neq \ica$, the distribution of values of $\hica$ 
can be quite complicated.  In Figure~\ref{fig:hilbvalues}, we plot these 
values for the matrix 
$$A = \left( \begin{array}{rrrr}
    1 & 3 & 2 & 4 \\
    1 & 4 & 5 & 2 \end{array} \right).$$
The boundary of $\cone(\A)$
is shown by dashed lines.  Notice that deep in the interior of the
cone, all of the values are one. Theorem~\ref{thm:hilbfunctionone}
proves this fact.

\begin{figure}
\begin{center}
\epsfig{file=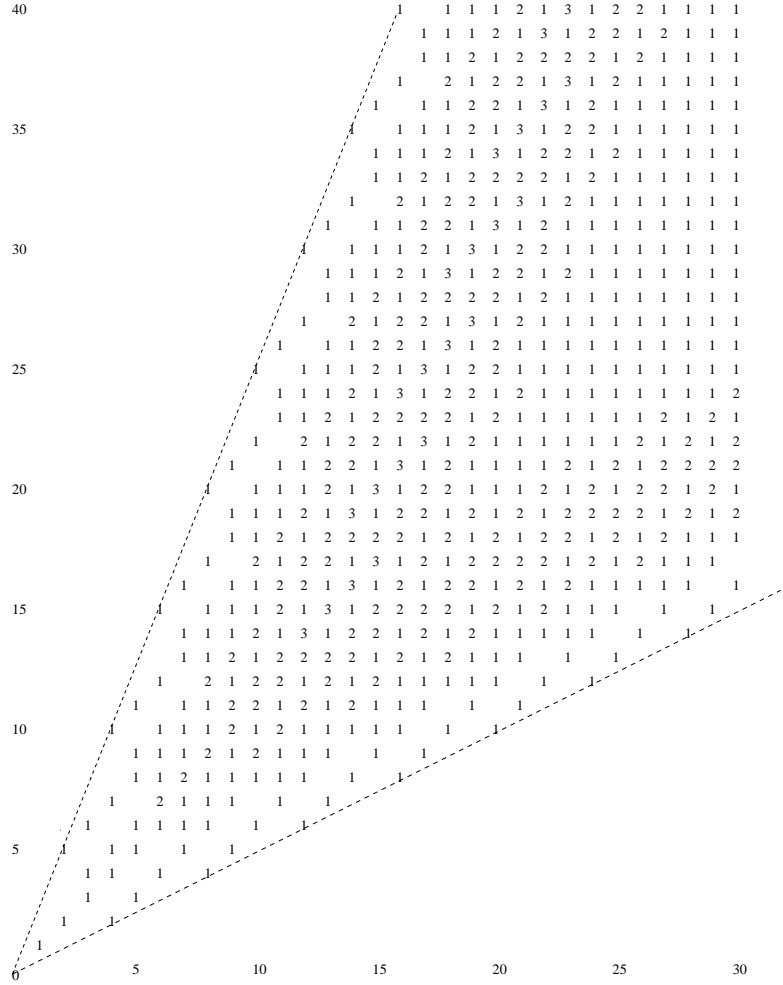,height=13.0cm} 
\end{center}
\caption{The distribution of values of $\hica$ for the matrix 
$A$ in Example~\ref{hf picture}.}
\label{fig:hilbvalues}
\end{figure}

\end{example}

\section{Monomial Initial Ideals of the Circuit Ideal}
\label{sec:initialideals}

Fix a generic weight vector $\omega \in \R^n$ such that
$\inomega(\ia)$ and $\inomega(\ica)$ are both monomial ideals. The
main result of this section is Theorem~\ref{thm:stdpairs} which
characterizes the associated primes of $\inomega(\ica)$ in terms of
certain polytopes defined from $\A$ and $\omega$ and their lattice
points. This theorem generalizes Theorem~2.5 in \cite{HT2} which gave
a complete characterization of the associated primes of
$\inomega(\ia)$ in terms of certain {\em lattice-point-free} polytopes
defined from $\A$ and $\omega$. Using Theorem~\ref{thm:stdpairs}, we
describe the similarities and differences between the associated
primes (standard pairs) of $\inomega(\ia)$ and $\inomega(\ica)$, and
give an answer to Problem~\ref{equality}
(Theorem~\ref{thm:obstruction}).

To begin, we prove that $\inomega(\ia)$ and $\inomega(\ica)$ have the
same radical. This result was stated in \cite{Ohs} without proof.
If $\u$ and $\v$ are vectors in $\Z^n$, then we say that $\u$ is
\emph{conformal} to $\v$ if $\supp(\u^+) \subseteq \supp(\v^+)$ and
$\supp(\u^-) \subseteq \supp(\v^-)$.
Lemma 4.10 in \cite{GBCP} states that every vector $\v \in \la$ can be 
written as a non-negative rational combination of $n-d$ circuits of $\A$ 
that are conformal to $\v$.


\begin{lemma}
\label{lem:circuitscaling}
If $\x^\alpha-\x^\beta\in \ia$ then there exists $s\in\N_{>0}$ such
that $\x^{s\alpha}-\x^{s\beta}\in \ica$.
\end{lemma}

\begin{proof}
  Suppose $\supp(\alpha) \cap \supp(\beta)=\emptyset$.  By \cite[Lemma
  4.10]{GBCP}, $\alpha-\beta = \sum_{i=1}^{t}q_i\cc_i$ where $\cc_i$
  are circuits of $\A$ conformal to $\alpha-\beta$ and $q_i\in\Q_{>
    0}$.  Clearing denominators and repeating the circuits in the sum
  if needed, there exists an $s\in\N_{>0}$ such that
  $s(\alpha-\beta)=\sum_{i=1}^{t'}\cc_i$.  Since the $\cc_i$ are
  conformal to $\alpha - \beta$, $s\alpha=\sum_{i=1}^{t'}{\cc_i}^+$
  and $s\beta=\sum_{i=1}^{t'}{\cc_i}^-$. Further, since
  $\x^{{\cc_i}^+}-\x^{{\cc_i}^-}\in \ica$ for each $i = 1, \ldots,
  t'$, $\x^{s\alpha}-\x^{s\beta}\in \ica$.  If $\supp(\alpha) \cap
  \supp(\beta) \not =\emptyset$ then by applying the above argument to
  $\x^{(\alpha-\beta)^+}-\x^{(\alpha-\beta)^-}$ we get
  $\x^{s(\alpha-\beta)^+}-\x^{s(\alpha-\beta)^-}\in \ica$. This
  implies that $\x^{s\alpha}-\x^{s\beta}
  =\x^{s(\min(\alpha,\beta))}(\x^{s(\alpha-\beta)^+}-
  \x^{s(\alpha-\beta)^-}) \in \ica$ where $\min(\alpha,\beta)$ is the
  component-wise minimum of $\alpha$ and $\beta$.
\end{proof}

\begin{proposition}
\label{prop:radicalsagree}
The radical ideals $\sqrt{\inomega(\ia)}$ and $\sqrt{\inomega(\ica)}$
coincide. 
\end{proposition}

\begin{proof}
  Since $\ia \supseteq \ica$, $\inomega(\ia) \supseteq \inomega(\ica)$
  and hence $\sqrt{\inomega(\ia)} \supseteq \sqrt{\inomega(\ica)}$.
  For the reverse inclusion it suffices to show that any monomial
  $\x^\u$ in $\sqrt{\inomega(\ia)} $ is in $\sqrt{\inomega(\ica)}$.
  If $\x^\u \in \sqrt{\inomega(\ia)} $ then $(\x^\u)^t\in
  \inomega(\ia)$ for some $t\in\N_{>0}$. Hence $\x^{t\u} =
  \inomega(\x^\alpha-\x^\beta)$ for some $\x^\alpha-\x^\beta\in \ia$.
  By Lemma \ref{lem:circuitscaling}, $\x^{s\alpha}-\x^{s\beta}\in
  \ica$ for some $s\in\N_{>0}$.  Thus,
  $\x^{ts\u}=\inomega(\x^{s\alpha}-\x^{s\beta})\in \inomega(\ica)$,
  and $\x^\u \in \sqrt{\inomega(\ica)}$.
\end{proof}

\begin{example} 
\label{ex:twistedcubic} 
Consider the matrix 
$$ A = \left( \begin{array}{rrrr} 
1 & 1 & 1 & 1\\
0 & 1 & 2 & 3
\end{array} \right).$$
Using the program \texttt{Gfan} \cite{Gfan} we find that both $\ia$
and $\ica$ have eight distinct monomial initial ideals.
Table~\ref{tab:radical} gives a representative weight vector $\omega$
for each pair of initial ideals and verifies
Proposition~\ref{prop:radicalsagree}.

\begin{table}

\begin{tabular}{|c|c|c|c|} \hline
 
$\omega$ & $\inomega(\ia)$ & $\inomega(\ica)$ 
& radical of both initial ideals \\ \hline

$(10,0,1,3)$ & $\langle ac, ad, bd \rangle$ 
& $\langle ac, a^2d, bd, ad^2 \rangle$ & $\langle ac, ad, bd \rangle$ \\

$(10,0,3,1)$ & $\langle ac, c^2, ad \rangle$ 
& $\langle ac, c^2, a^2d, abd, ad^2 \rangle$  & $\langle c, ad \rangle$ \\ 

$(3,1,10,0)$ & $\langle ac, bc, c^2, a^2d \rangle$ & $\langle ac,
b^2c, c^2, a^2d, bcd \rangle$ & ''\\

$(1,3,10,0)$ & $\langle b^3, ac, bc, c^2 \rangle$ 
& $\langle b^3, ac, b^2c, c^2, bcd \rangle$ & $\langle b, c\rangle$ \\

$(1,5,3,0)$ & $\langle b^2, bc, c^2 \rangle$ 
& $\langle b^2, abc, c^2, bcd \rangle$ & '' \\

$(0,10,3,1)$ & $\langle b^2, bc, c^3, bd \rangle$ 
& $\langle b^2, abc, bc^2, c^3, bd \rangle$ & '' \\ 

$(1,3,0,10)$ & $\langle b^2, ad, bd \rangle$ 
& $\langle b^2, a^2d, bd, acd, ad^2 \rangle$ & $\langle b, ad \rangle$ \\

$(3,10,0,1)$ & $\langle b^2, bc, bd, ad^2 \rangle$ 
& $\langle b^2, abc, bc^2, bd, ad^2 \rangle$ & '' \\
\hline
\end{tabular}

\vspace{.2cm}

\caption{Comparison of initial ideals of $\ia$ and $\ica$ from
  Example~\ref{ex:twistedcubic}.} 
\label{tab:radical}
\end{table}
\end{example}

\vspace{.2cm}

\begin{definition}~\cite{GBCP}
\begin{enumerate}
\item The \emph{regular triangulation} of $\A$ with respect to
  $\omega$ is the simplicial complex $\Delta_{\omega}$ on the vertex
  set $[n] = \{1, \ldots, n\}$ such that $\{i_1,\ldots,i_r\}
  \subseteq [n]$ is a face of $\Delta_{\omega}$ if and only if there
  exists a vector $\cc \in \R^d$ such that $\a_j \cdot \cc = \omega_j$
  if $j \in \{i_1,\ldots,i_r\}$ and $\a_j \cdot \cc < \omega_j$ if $j
  \notin \{i_1,\ldots,i_r\}$.
  
\item The \emph{Stanley-Reisner} ideal of a simplicial complex
  $\Delta$ on $[n]$ is the ideal in $\k[\x]$ generated by the
  monomials $\x_{\sigma} := \prod_{i \in \sigma} x_i$ for each minimal
  nonface $\sigma$ of $\Delta$.
\end{enumerate}
\end{definition}

Theorem 8.3 in \cite{GBCP} states that $\sqrt{\inomega(\ia)}$ is the
Stanley-Reisner ideal of the regular triangulation $\Delta_{\omega}$
of $\A$.  For a set $\sigma \subseteq [n]$ define $P_{\sigma} :=
\langle x_j \, : \, j \not \in \sigma \rangle \subset \k[\x]$.  Note
that $P_{\sigma}$ is a monomial prime ideal such that
$\k[\x]/P_\sigma$ has Krull dimension $|\sigma|$.

\begin{corollary} \label{cor:initialassprimes}
\begin{enumerate}
\item All the associated primes of $\inomega(\ica)$ are
  monomial ideals of the form $P_{\sigma}$ where $\sigma$ is a face of
  the simplicial complex $\Delta_{\omega}$.
\item The prime $P_{\sigma}$ is a minimal prime of
  $\inomega(\ica)$ if and only if $\sigma$ is a maximal
  face of $\Delta_{\omega}$. 
\item The ideal $\inomega(\ica)$ is equi-dimensional.
\end{enumerate}
\end{corollary}

\begin{proof} 
  If $I$ is the Stanley-Reisner ideal of a simplicial complex $\Delta$
  on $[n]$, then $I$ has the irredundant prime decomposition $I =
  \cap_{\sigma \in \max(\Delta)} P_{\sigma}$ where $\max(\Delta)$ is
  the set of maximal faces of $\Delta$~\cite[Chapter 8]{GBCP}.  Thus the 
  minimal primes of $\inomega(\ica)$, which equal the minimal primes of
  $\sqrt{\inomega(\ica)} = \sqrt{\inomega(\ia)}$, are the primes
  $P_{\sigma}$ as $\sigma$ varies in $\max(\Delta_{\omega})$, proving
  (2).  Since $\Delta_{\omega}$ is a pure simplicial complex, we get
  (3).  If $P_{\tau}$ is an embedded prime of $\inomega(\ica)$, then
  $\tau \subset \sigma$ for some $\sigma \in \max(\Delta_{\omega})$.
  This implies that $\tau$ is a lower dimensional face of
  $\Delta_{\omega}$, proving (1).
\end{proof}

If $\tau$ is a lower dimensional face of $\Delta_{\omega}$, $P_{\tau}$
may or may not be an embedded prime of $\inomega(\ica)$.
Theorem~\ref{thm:stdpairs} characterizes the lower dimensional faces
of $\Delta_{\omega}$ that index embedded primes of $\inomega(\ica)$.

\begin{remark} \label{rem:arbitraryspanning} 
  If $\B$ is an arbitrary spanning set of $\la$ and $J_{\B}$ as in
  Proposition~\ref{prop:saturation}, then it need not be that
  $\sqrt{\inomega(J_{\B})}$ is the Stanley-Reisner ideal of any
  regular triangulation of $\A$. In Example~\ref{ex:twistedcubic}, the
  set $\B = \{ (1,-2,1,0), (2,-3,0,1) \}$ spans the lattice $\la$ and
  $J_{\B} = \langle ac-b^2, a^2d-b^3 \rangle$. The grevlex initial
  ideal of $J_{\B}$ with $a \succ b \succ c \succ d$ is $\langle b^2,
  abc, a^2c^2 \rangle$ whose radical is $\langle b, ac \rangle$.  This
  ideal is not listed in the last column of Table~\ref{tab:radical}.
\end{remark}

We now establish the necessary definitions and lemmas
for Theorem~\ref{thm:stdpairs}. The associated primes of a monomial
ideal $M$ can be studied via a combinatorial construction introduced
in \cite{STV} called the {\em standard pair decomposition} of $M$.

\begin{definition} \label{def:stdpairs}
  Let $M \subset \k[\x]$ be a
  monomial ideal, $\x^\u$ a standard monomial of $M$ and
  $\sigma \subseteq [n]$. Then $(\x^\u, \sigma)$ is an 
\emph{admissible pair} of $M$ if:
\begin{enumerate}
\item $\supp(\u) \cap \sigma = \emptyset$,
\item all monomials in $\x^\u \cdot {\bf k}[x_j \, : \, j
  \in \sigma]$ are standard monomials of $M$.
\end{enumerate}
An admissible pair $(\x^\u, \sigma)$ of $M$ is called a 
{\em standard pair} of $M$ if there does not exist another 
admissible pair $(\x^\v, \tau)$ such that 
$\v \leq \u$ and $\supp(\u - \v) \cup
\sigma \subseteq  \tau$. 
\end{definition}

The (unique) decomposition of the
standard monomials of $M$ given by its standard pairs is the
\emph{standard pair decomposition} of $M$. Let $\ass(I)$ denote the
set of associated primes of an ideal $I$. Since $M$ is a monomial
ideal, all elements of $\ass(M)$ have the form $P_{\sigma}$ for some
$\sigma \subseteq [n]$. Standard pairs of $M$ are related to $\ass(M)$
as follows.

\begin{proposition} \cite{STV} \label{prop:stdpairsandmonomials}
\begin{enumerate}
\item $P_{\sigma} \in \ass(M)$ if and only if
  $M$ has a standard pair of the form $(\cdot, \sigma)$.
\item $P_{\sigma}$ is a minimal prime of $M$ if and only if 
$(1, \sigma)$ is a standard pair of $M$.
\end{enumerate}
\end{proposition}

\begin{figure}[ht]
\input{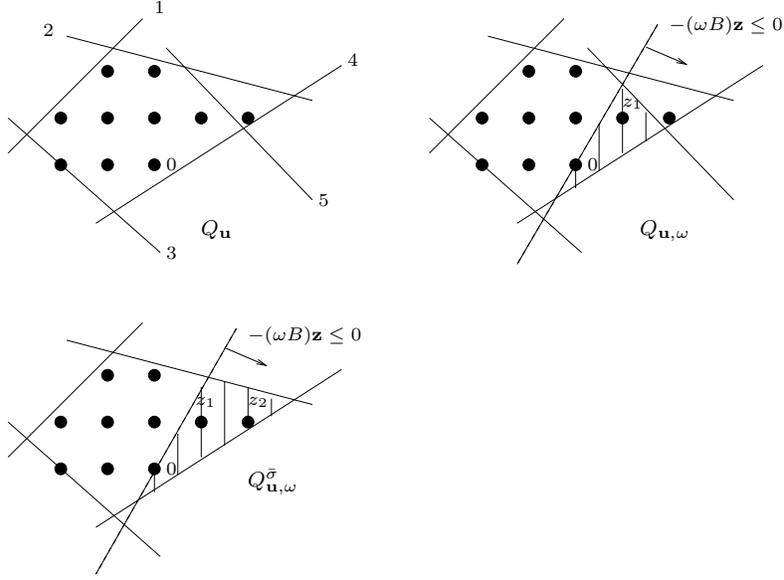}
\caption{The polytopes $Q_{\u}$, $Q_{\u,\omega}$ and
  $Q_{\u,\omega}^{\bar \sigma}$. \label{fig:qpolys}}
\end{figure}

We now define the polytopes needed in Theorem~\ref{thm:stdpairs}.  Fix
a matrix $B \in \Z^{n \times (n-d)}$ whose columns form a basis for
the lattice $\la$.  Such a $B$ is called a {\em Gale dual} of $A$. In
particular, the columns of $B$ span the kernel of $A$ as an
$\R$-vector space.  For $\u \in \N^n$ let
$$Q_{\u} := \{ \z \in \R^{n-d} \, : \, B\z \leq \u \}.$$
Recall that
by assumption, $P_{\b} = \{\x \in \R^n_{\geq 0} \, : \, A\x = \b \}$
is a polytope for all $\b \in \na$. The polyhedron $Q_\u$ is the image
of $P_{A\u}$ under the isomorphism $$\phi_\u \, : \, \{\x \in \R^n \,
: \, A\x = A \u \} \rightarrow \R^{n-d} \,\,\textup{such that} \,\, \x
\mapsto \z \,\,\textup{where} \,\,\u - B\z = \x.$$
For each $\x \in
\R^n$ such that $A\x = A\u$, $\u-\x = B \z$ for some $\z \in \R^{n-d}$
since $\u-\x \in \ker(A) = \{ B\z \, :\, \z \in \R^{n-d} \}$.
Further, this $\z$ is unique since the columns of $B$ are linearly
independent. The vector $\u$ maps to ${\bf 0}$ under $\phi_{\u}$ and
hence ${\bf 0} \in Q_{\u}$.

Next, define $$Q_{\u, \omega} := Q_{\u} \cap \{ \z \in \R^{n-d} \, :
\, (-\omega B) \z \leq 0 \}, $$
the subpolytope of $Q_{\u}$ created by
adding one new inequality depending on $\omega$.  For $\sigma$ a face
of $\Delta_{\omega}$, further define
$$Q_{\u, \omega}^{\bar \sigma} := \{ \z \in \R^{n-d} \, : \, (B\z \leq
\u)^{\bar \sigma}, (-\omega B) \z \leq 0 \}$$
where $(B\z \leq
\u)^{\bar \sigma}$ denotes the subsystem of inequalities indexed by
$\bar \sigma$ in $B \z \leq \u$. Theorem~1 in \cite{HT01} proves that
$Q_{\u,\omega}^{\bar \sigma}$ is a polytope. It is a relaxation of
$Q_{\u, \omega}$.  Figure~\ref{fig:qpolys} shows pictures of the
polytopes $Q_{\u}$, $Q_{\u,\omega}$ and $Q_{\u,\omega}^{\bar \sigma}$.
The inequalities in $B\z \leq \u$ are numbered $1,\ldots,5$ and
$Q_{\u,\omega}^{\bar \sigma}$ is drawn for $\sigma = \{5\}$.

For a non-zero lattice point $\z \in Q_{\u, \omega}^{\bar \sigma}$, set 
$\m_\z := (\u - B\z)^{-}$. Let $B_i$ denote the $i$-th row of $B$.  

\begin{remark}
\begin{enumerate}
\item The $i$-th component $(\m_{\z})_i > 0$ if and only if $\z$ violates 
the $i$-th inequality $B_i\z \leq u_i$ among the inequalities 
$B\z \leq \u$ defining $Q_{\u, \omega}$. 
\item Since every $\z \in Q_{\u, \omega}^{\bar \sigma}$
satisfies $B_i \z \leq u_i$ for $i \not \in \sigma$, the support of
$\m_\z$ is contained in $\sigma$.  
\item The vector $\m_\z$ is the component-wise smallest vector $\m$ in 
$\N^n$ with support in $\sigma$ such that $\z \in Q_{\u + \m, \omega}$.  
\item By the definition of $\m_\z$, $\u+\m_{\z} - B\z \in \N^n$. 
\end{enumerate}
\end{remark}

For instance, in Figure~\ref{fig:qpolys}, $\z_1 \in Q_{\u,\omega}$ and hence
$\m_{\z_1} = {\bf 0}$, while $\z_2$ violates the inequality 
$B_5 \z \leq u_5$ defining $Q_{\u,\omega}$ and hence $\m_{\z_2}$ 
has a positive fifth component.

Theorem~\ref{thm:stdpairs} will generalize the following theorem for
toric ideals.

\vspace{.2cm}

\begin{theorem}~\cite{HT2} \label{thm:toriccase}
  Assume $\u \in \N^n$ and $\sigma \in \Delta_{\omega}$ such that
  $\supp(\u) \cap \sigma = \emptyset$. Then $(\x^\u, \sigma)$ is a
  standard pair of $\inomega(\ia)$ if and only if the following two
  conditions hold.
\begin{enumerate}
\item There are no non-zero lattice points in $Q_{\u,
    \omega}^{\bar \sigma}$.
\item For every $i \in \bar{\sigma}$ there is a non-zero lattice point
  in $Q_{\u,\omega}^{\overline{\sigma \cup \{i\}}}$.
\end{enumerate}
\end{theorem}

Theorem~\ref{thm:stdpairs} is analogous, but involves an algebraic
component rather than being purely polyhedral.  Recall that
$\x_{\sigma} = \prod_{i \in \sigma} x_i$.

\begin{theorem} \label{thm:stdpairs}
  Assume $\u \in \N^n$ and $\sigma \in \Delta_{\omega}$ such that
  $\supp(\u) \cap \sigma = \emptyset$. Then $(\x^{\u}, \sigma)$ is a
  standard pair of $\inomega(\ica)$ if and only if the following two
  conditions hold.
\begin{enumerate}
\item For each non-zero lattice point $\z$ in $Q_{\u, \omega}^{\bar
    \sigma}$, $\x^{\u+\m_{\z}} - \x^{\u+\m_{\z} - B\z} \not \in (\ica
  \, : \, \x_{\sigma}^{\infty})$.
\item For each $i \in \bar{\sigma}$, there exists some non-zero
  lattice point $\z \in Q_{\u, \omega}^{\overline {\sigma \cup
      \{i\}}}$ such that $\x^{\u+{\bf m}_{\z}} - \x^{\u+\m_{\z} - B\z}
  \in (\ica \, : \, \x_{\sigma \cup \{i\}}^{\infty})$.
\end{enumerate}
\end{theorem}

\begin{remark} \label{rem:Jstdpairs}
  Let $J$ be any ideal such that $\ica \subseteq J \subseteq \ia$.
  Corollary~\ref{cor:initialassprimes} and Theorem~\ref{thm:stdpairs}
  apply to $\inomega(J)$ by simply replacing $\ica$ by $J$ everywhere
  in the statements and proofs.
\end{remark}

We first use Theorem~\ref{thm:stdpairs} to reprove
Theorem~\ref{thm:toriccase}.

\smallskip\noindent {\it Proof of Theorem~\ref{thm:toriccase}}: Since
$\ia$ is prime and monomial free, $(\ia : \, \x_{\tau}^{\infty}) =
\ia$ for all $\tau \subseteq [n]$. Thus if $\z$ is a non-zero lattice
point in $Q_{\u, \omega}^{\bar \sigma}$, then $\x^{\u+\m_\z} -
\x^{\u+\m_\z - B\z} \in \ia = (\ia \, : \, \x_{\sigma}^{\infty})$.
Hence, Theorem~\ref{thm:stdpairs} (1) holds if and only if there are
no non-zero lattice points in $Q_{\u, \omega}^{\bar \sigma}$.
Similarly, Theorem~\ref{thm:stdpairs} (2) holds in the toric situation
if and only if for every $i \in \bar{\sigma}$ there is a non-zero
lattice point in $Q_{\u,\omega}^{\overline{\sigma \cup \{i\}}}$.
\endproof

\smallskip\noindent {\it Proof of Theorem~\ref{thm:stdpairs}}:
$(\Rightarrow)$: Suppose $(\x^\u, \sigma)$ is a standard pair of
$\inomega(\ica)$. Then $\sigma \in \Delta_{\omega}$ and $\supp(\u)
\cap \sigma = \emptyset$. 
Suppose $\z$ is a non-zero lattice point in $Q_{\u, \omega}^{\bar \sigma}$.  
Then $- (\omega B)\z \leq 0$, and because $\omega$ is generic, we may
assume $- (\omega B)\z < 0$. For any $\m \in \N^n$ with support
contained in $\sigma$, $\x^{\u + \m}$ is a standard monomial of
$\inomega(\ica)$ since $(\x^\u, \sigma)$ is a standard pair. If
further, $\m \geq \m_{\z} = (\u - B\z)^{-}$, then 
$\u +\m - B\z \in \N^n$ and $A(\u + \m - B\z) = A(\u + \m)$ since $AB
= 0$.  Also, $ \omega \cdot (\u + \m - B\z) = \omega \cdot (\u + \m)
- (\omega B)\z < \omega \cdot (\u + \m)$ since $- (\omega B)\z < 0$.
Therefore, by Lemma~\ref{lem:cantconnect}, $\x^{\u + \m} - \x^{\u + \m
  - B\z} \not \in \ica$.  In particular, $\x^{\u + \m_\z} - \x^{\u +
  \m_\z - B\z} \not \in \ica$ and for all $\m' \in \N^n$ with support
in $\sigma$, $\x^{\bf m'}(\x^{\u + \m_\z} - \x^{\u + \m_\z - B\z})
\not \in \ica$.  Rewriting, this gives $\x^{\u + \m_\z} - \x^{\u +
  \m_\z - B\z} \not \in (\ica \, : \, \x_{\sigma}^{\infty})$ and (1)
holds.

Suppose $i \not \in \sigma$. Then there exists some $\m \in \N^n$ with
$\supp(\m) \subseteq \sigma$ and $p > 0$ such that $\x^{\u+\m}x_i^p
\in \inomega(\ica)$. Let $\q$ be the unique sink in the same component
of $\F_{\omega}^{\C}(A(\u+\m+p\e_i))$ as $\u + \m + p\e_i$.  Note that
$\q \neq \u + \m + p\e_i$ since $\x^\q \not \in \inomega(\ica)$. Let
$\z \in \Z^{n-d}$ be such that $\q= \u + \m + p\e_i - B\z$. Then $\u +
\m + p\e_i$ maps to ${\bf 0}$ and $\q$ maps to $\z \neq {\bf 0}$ in
$Q_{\u + \m + p\e_i}$ under the map $\phi_{\u + \m + p\e_i}$.  
Since $\omega \cdot \q= \omega \cdot (\u + \m
+ p\e_i - B\z) < \omega \cdot (\u + \m + p\e_i)$, we see that
$-(\omega B)\z < 0$.  Therefore, $\z$ is a lattice point in $Q_{\u +
  \m + p\e_i, \omega}$ and hence in $Q_{\u,\omega}^{\overline{\sigma
    \cup \{i\}}}$ obtained by throwing away the inequalities of $B\z
\leq \u$ indexed by $\sigma \cup \{i\}$ from $Q_{\u + \m + p\e_i,
  \omega}$.  This is because $\supp(\m+p\e_i) \subseteq \sigma \cup \{i\}$.
By definition, $\m_{\z} \leq \m + p\e_i$
since $\m_{\z}$ is the component-wise smallest vector $\m'$ with
support in $\sigma \cup \{i\}$ such that $\z \in
Q_{\u+\m',\omega}^{\overline{\sigma \cup \{i\}}}$ and we know that 
$\z \in Q_{\u + \m + p\e_i}$.  Since $\q= {\u +
  \m + p\e_i-B\z}$ lies in the same component of
$\F_{\omega}^{\C}(A(\u+\m+p\e_i))$ as $\u + \m + p\e_i$, by
Lemma~\ref{lem:samecomponent},
  $$ \x^{\u + \m + p\e_i} - \x^{\u + \m + p\e_i-B\z} 
  = \x^{\m + p\e_i - \m_{\z}}(\x^{\u + \m_{\z}} - 
  \x^{\u + \m_{\z}-B\z}) \in \ica.$$
  This implies that $\x^{\u + \m_{\z}} - \x^{\u + \m_{\z}-B\z} 
  \in (\ica \, :\, \x_{\sigma \cup \{i\}}^{\infty})$ and (2) holds.
  
  $(\Leftarrow)$: Suppose (1) and (2) hold for some $\sigma \in
  \Delta_{\omega}$ and some $\u \in \N^n$ with support in $\bar
  \sigma$.  We first show that $\x^{\u+\m}$ is a standard monomial of
  $\inomega(\ica)$ where $\m \in \N^n$ is an arbitrary vector with
  $\supp(\m) \subseteq \sigma$.  
  Suppose $\z$ is a non-zero lattice point in $Q_{\u+\m, \omega}$. Then
  $\z$ is also a non-zero lattice point in the relaxation
  $Q_{\u+\m,\omega}^{\bar \sigma} = Q_{\u,\omega}^{\bar \sigma}$.
  Compute $\m_\z$ for this $\u$ and $\z$. Since $\z \in Q_{\u+\m,
    \omega}$, $\m_{\z} \leq \m$. By (1), $(\x^{\u + \m_\z} - \x^{\u +
    \m_\z - B\z}) \not \in (\ica \, : \, \x_{\sigma}^{\infty})$ which
  implies that
  $$\x^{\m-\m_{\z}}(\x^{\u + \m_\z} - {\bf x}^{\u + \m_\z - B\z}) =
  \x^{\u + \m} - \x^{\u + \m - B\z} \not \in \ica.$$
  Thus for each $\z \neq {\bf 0}$ in $Q_{\u+\m, \omega}$, the vector 
  $\u+\m-B\z$ does not lie in the same component as $\u+\m$. This 
  implies that $\omega \cdot \v > \omega \cdot (\u+\m)$ 
  for all $\v$ in the same component
  as $\u+\m$. By Lemma~\ref{lem:cantconnect}, $\x^{\u+\m}$ is a
  standard monomial of $\inomega(\ica)$.  Since $\supp(\u) \cap \sigma
  = \emptyset$ and $\m$ is an arbitrary vector with support contained
  in $\sigma$, we conclude that $(\x^\u, \sigma)$ is an admissible
  pair of $\inomega(\ica)$.
  
  To show that $(\x^\u, \sigma)$ is a standard pair, we need
  to argue that the monomials of this pair are not properly contained
  in any other standard pair $(\x^{\bf u'}, \sigma')$ of
  $\inomega(\ica)$.  Suppose there is such a standard pair.  We
  first argue that $\sigma = \sigma'$. By (2), if $i \not \in \sigma$
  then there exists some non-zero lattice point $\z$ in 
  $Q_{\u, \omega}^{\overline{\sigma \cup \{i\}}}$ such that
  $$
  \x^{\u+\m_\z} - \x^{\u+\m_\z - B\z} \in (\ica \, : \, \x_{\sigma
    \cup \{i\}}^{\infty}).$$
  This implies that there exists some $p
  \in \N$ and $\m \in \N^n$ with support in $\sigma$ such that
  $x_i^p\x^\m(\x^{\u+\m_\z} - \x^{\u+\m_\z - B\z}) \in \ica$.  Since
  $(-\omega B)\z < 0$, $x_i^p\x^\m(\x^{\u+\m_\z})$ is the leading term
  of $x_i^p\x^\m(\x^{\u+\m_\z} - \x^{\u+\m_\z - B\z}) \in \ica$ and
  hence is in $\inomega(\ica)$.  This construction shows that not all
  monomials of the form $\x^\u \x^\q$ where the support of ${\bf q}$
  is contained in $\sigma \cup \{i\}$ are standard monomials of
  $\inomega(\ica)$ and hence $(\x^\u, \sigma)$ is not contained in any
  admissible pair $(\x^\u, \sigma')$ with $\sigma \subsetneq \sigma'$.
  To finish the argument, suppose $(\x^\u, \sigma)$ is contained in a
  standard pair of form $(\x^{\bf u'}, \sigma)$.  Then $\u = \m\u'$
  for some $\m$ whose support is contained in $\sigma$.  However,
  because $(\x^\u, \sigma)$ is a standard pair, the support of $\u$
  must also be disjoint from $\sigma$.  Thus $\m = 1$ and so $\u =
  \u'$.  \endproof

We now apply Theorems~\ref{thm:stdpairs} and~\ref{thm:toriccase} to
study the difference between the two monomial ideals 
$\inomega(\ia)$ and $\inomega(\ica)$.  This difference will be the
key to our study of the associated primes of $\ica$ itself in 
Section~\ref{sec:assprimes}.

\begin{definition}\label{def:cspandtsp}
  A \emph{circuit-specific standard pair} (CSP) is a standard pair of 
  $\inomega (\ica)$ that is not also a standard pair of $\inomega (\ia)$.
\end{definition}

\begin{corollary} \label{cor:CSPchar} 
  Assume $\u \in \N^n$ and $\sigma \in \Delta_{\omega}$ such that
  $\supp(\u) \cap \sigma = \emptyset$.  Then $(\x^\u, \sigma)$ is a
  CSP if and only if the two conditions of Theorem~\ref{thm:stdpairs}
  hold and there exists at least one non-zero lattice point $\z \in
  Q_{\u, \omega}^{\bar \sigma}$.
\end{corollary}

\begin{proof}
  If the two conditions of Theorem~\ref{thm:stdpairs} hold then
  $(\x^\u,\sigma)$ is a standard pair of $\inomega (\ica)$ and if there
  is a non-zero lattice point $\z \in Q_{\u, \omega}^{\bar
    \sigma}$, then by Theorem~\ref{thm:toriccase}, $(\x^\u,\sigma)$
  is not a standard pair of $\inomega (\ia)$. Thus it is a CSP.
  Conversely, if $(\x^\u,\sigma)$ is a CSP then the two conditions of
  Theorem~\ref{thm:stdpairs} hold.  Suppose there is no nonzero lattice point 
  $\z \in Q_{\u, \omega}^{\bar \sigma}$.  Then condition (1) of
  Theorem~\ref{thm:toriccase} is true. But since condition (2) of
  Theorem~\ref{thm:stdpairs} holds for this CSP, there is a non-zero
  lattice point in $Q_{\u,\omega}^{\overline{\sigma \cup \{i\}}}$
  for each $i \not \in \sigma$, which is condition (2) of
  Theorem~\ref{thm:toriccase}. This implies that $(\x^\u,\sigma)$ is a
  standard pair of $\inomega(\ia)$, contradicting that it is a CSP.
\end{proof}

\begin{example} Consider the matrix $A$ and weight vector $\omega$
  given below:
$$A = \left( 
\begin{array}{ccccc} 
1 & 1 & 1 & 1 & 1 \\
0 & 3 & 4 & 5 & 6 \\
0 & 0 & 7 & 8 & 9 
\end{array}
\right), \,\,\,\,\, \omega = (1000,100,10,1,0).$$
The matrix $$ B = \left( \begin{array}{rr} 
2& 4\\
-1 & -2\\
-5 & -9\\
1 &  0\\
3 &  7
 \end{array} \right)$$ 
is a Gale dual of $A$ and $\omega B = (1851, 3710)$. The circuit ideal 
$$\ica = \langle b^2c^9 - a^4e^7 ,bd^9 - a^2e^8 ,bc^8 - a^2d^7 ,ce -
d^2 \rangle$$
and its initial ideal $$
\inomega(\ica) = \langle
a^2d^7, ce, a^4d^6e^4, a^4d^4e^5, a^2d^6e^5, a^4d^2e^6, a^2d^4e^6,
a^4e^7, a^2d^2e^7, a^2e^8 \rangle $$
has $58$ standard pairs. These
ideals and standard pairs were computed using \texttt{Macaulay
  2}~\cite{M2}. Consider the standard pair $(d^4e^3,
\{1,2\})$ for which $\u = (0,0,0,4,3)$ and $\sigma = \{1,2\}$. The
monomial $d^4e^3$ is a standard monomial for the toric initial ideal
$\inomega(\ia)$ as well and $Q_{\u,\omega} \cap \Z^2 = \{0\}$.
However, the polytope
$$
Q_{\u, \omega}^{\bar \sigma} = 
\left\{ \z \in \Z^2 \, : \, \left(
\begin{array}{rr}
-5 & -9\\
1 & 0\\
3 & 7
\end{array} \right)
\z \leq 
\left(
\begin{array}{r}
0 \\ 4 \\ 3
\end{array} \right), \,\, 
1851 z_1 + 3710 z_2 \geq 0 \right\}$$
contains two more lattice
points: $(1,0)$ and $(3,-1)$.  Thus $(\x^\u, \sigma)$ is not a
standard pair of $\inomega(\ia)$, so it is a CSP.  Both points have
$\m_\z = (2,0,0,0,0)$.  For $(1,0)$, $\x^{\u+\m_\z} - \x^{\u+\m_\z -
  B\z} = a^2d^4e^3-bc^5d^3$ is not in $(\ica \, : \,(ab)^{\infty})$
but does lie in $(\ica \, : \, (aby)^{\infty})$ for each $y \in
\{c,d,e\}$. Similarly, for $(3,-1)$, $\x^{\u+\m_\z} - \x^{\u+\m_\z -
  B\z} = a^2d^4e^3-bc^6de$ which also does not lie in $(\ica \, : \,
(ab)^{\infty})$ but does lie in $(\ica \, : \, (aby)^{\infty})$ for
each $y \in \{c,d,e\}$. 
\end{example} 

We now use CSPs to give a precise description of the set $\mathcal B
:= \{ \b \in \na \,:\, \hica(\b) > 1 \}$. This description gives a 
new proof of the following theorem alluded to in 
Section~\ref{sec:optimization} (cf. Figure~\ref{fig:hilbvalues}).  The
theorem also follows from~\cite[Corollary~5.3]{DES}.

\begin{theorem}
\label{thm:hilbfunctionone}
For all $\b \in \na$ sufficiently far from the boundary of $\cone(\A)$,
$\hica(\b) = 1$ and hence the graphs $\F^\C(\b)$ and
$\F_\omega^{\C}(\b)$ are connected.  
\end{theorem}

Recall that $\b$ lies in $\B$ if and only if for a generic $\omega$,
$\inomega(\ica)$ has more than one standard monomial of degree $\b$.
That is, $\mathcal B = \{ A \u \,:\, \x^\u \in \inomega(\ia) \backslash
\inomega(\ica) \}$.  Since all standard monomials of degree $\b$ other 
than the toric standard monomial lie on CSPs of $\inomega(\ica)$, it
follows that $\mathcal B$ is contained in the union of the images in $\na$ 
of the CSPs of $\inomega(\ica)$ under the map 
$\pi \,:\, \N^n \rightarrow \na$, $\u \mapsto A\u$.  

\begin{lemma}
\label{lem:nointeriorCSP}
If $(\x^\u, \sigma)$ is a CSP of $\inomega(\ica)$, then the set
$\{\a_i: i \in \sigma \}$ is contained in a facet of $\cone(\A)$. 
\end{lemma}

\begin{proof}
  Since $(\x^\u, \sigma)$ is a standard pair of $\inomega(\ica)$, by
  Corollary~\ref{cor:initialassprimes}, $\sigma$ is a face of the
  regular triangulation $\Delta_\omega$. Suppose $\cone(\A_{\sigma})$
  intersects the interior of $\cone(\A)$. Choose a monomial
  $\x^\alpha$ on the CSP $(\x^\u, \sigma)$ such that $\x^\alpha \in
  \inomega(\ia)$.  Let $\x^\beta$ be the standard monomial of
  $\inomega(\ia)$ of degree $A\alpha$. Then $\x^\alpha - \x^\beta \in
  \ia$ with leading term $\x^\alpha$. Since $\x_\sigma^\m \x^\alpha \not
  \in \inomega(\ica)$ for any $\m$, the binomial $\x_\sigma^\m (\x^\alpha
  - \x^\beta) \not \in \ica$ for any $\m$ since its leading term
  $\x_\sigma^\m \x^\alpha$ would then have to be in $\inomega(\ica)$.
  This implies that $(\ica \,:\, \x_\sigma^\infty) \neq \ia$. On the
  other hand, since every embedded prime of $\ica$ is of the form
  $P_\tau + \ia$ where $\tau$ indexes some proper face of
  $\cone(\A)$ (see Proposition~\ref{prop:primarydecomposition}), the
  monomial $\x_\sigma$ lies in each of these embedded primes since
  $\sigma$ is not contained in any proper face of $\cone(\A)$. This
  implies that for $\m$ large enough, $\x_\sigma^\m$ lies in every
  primary component of $\ica$ except $\ia$, which in turn
  implies that $(\ica \,:\, \x_\sigma^\infty) = \ia$, a contradiction.
\end{proof}

\noindent {\em Proof of Theorem~\ref{thm:hilbfunctionone}}:
  By Lemma~\ref{lem:nointeriorCSP}, if $(\x^\u, \sigma)$ is a CSP of
  $\inomega(\ica)$, then $A\u + \N \A_\sigma$, its image under $\pi$
  in $\na$, is contained in some hyperplane parallel to a facet of
  $\cone(\A)$.  Since there are only finitely many CSPs of
  $\inomega(\ica)$, $\mathcal B$ is contained in finitely many
  hyperplanes parallel to the facets of $\cone(\A)$. This implies that
  the maximum distance of a point in $\mathcal B$ from the boundary of 
  $\cone(\A)$ is bounded which proves the theorem. \hfill$\square$\medskip
  
  Using Lemma~\ref{lem:nointeriorCSP} we can also prove that even if
  $\inomega(\ia) \neq \inomega(\ica)$, they share a significant number
  of standard pairs. Applications of this result will be developed in
  Section~\ref{sec:fans}.  By Corollary~\ref{cor:initialassprimes},
  $\inomega(\ia)$ and $\inomega(\ica)$ have the same minimal primes
  $P_{\sigma}$ as $\sigma$ varies over $\max(\Delta_{\omega})$.  For
  such primes, $|\sigma| = d$ and $\cone(\A_\sigma)$ is
  $d$-dimensional.

\begin{definition} \label{def: fat stdpairs}
  A standard pair $(\x^{\u},\sigma)$ of $\inomega(\ia)$ or
  $\inomega(\ica)$ is said to be {\em fat} if $P_{\sigma}$
  is a minimal prime of these initial ideals or equivalently if 
  $|\sigma| = d$.
\end{definition}

 We now prove that $\inomega(\ia)$ and $\inomega(\ica)$ have the same
 fat standard pairs, strengthening
 Proposition~\ref{prop:radicalsagree}.  An alternate proof relying on
 enumeration is included in Section~\ref{sec:fans}.

\begin{theorem} \label{thm:fatstdpairsagree} 
  The monomial ideals $\inomega(\ia)$ and $\inomega(\ica)$ have the 
  same fat standard pairs.
\end{theorem}

\begin{proof}
  Suppose $(\x^\u, \sigma)$ is a fat standard pair of $\inomega(\ia)$.
  Since $\inomega(\ica) \subseteq \inomega(\ia)$, $(\x^\u, \sigma)$ is
  an admissible pair of $\inomega(\ica)$. If it is not a standard pair
  of $\inomega(\ica)$, there exists some standard pair $(\x^\v, \tau)$
  of $\inomega(\ica)$ such that $\supp(\u-\v) \cup \sigma \subseteq
  \tau$. However, since $\sigma \in \max(\Delta_\omega)$ and $\tau \in
  \Delta_\omega$, $\sigma = \tau$. This implies that $\supp(\u-\v)
  \subseteq \sigma$ which in turn implies that $\u = \v$ since
  $\supp(\u) \cap \sigma = \supp(\v) \cap \sigma = \emptyset$.
  
  Suppose $(\x^\u, \sigma)$ is a fat standard pair of
  $\inomega(\ica)$. Then by Proposition~\ref{prop:stdpairsandmonomials} 
  and Corollary~\ref{cor:initialassprimes}, $\sigma$ is a maximal 
  face of $\Delta_\omega$. Then by Lemma~\ref{lem:nointeriorCSP}, 
  $(\x^\u, \sigma)$ is not a CSP of $\inomega(\ica)$, so 
  $(\x^\u, \sigma)$ is a standard pair of $\inomega(\ia)$.
\end{proof}

We conclude this section with an answer to Problem~~\ref{equality}. 

\begin{definition} \label{def:CSPpolytope}
  A polytope $Q_{\u, \omega}^{\bar \sigma}$ corresponding to a CSP 
  $(\x^\u, \sigma)$ of $\inomega(\ica)$ is called a 
  {\em CSP polytope} of $\A$.
\end{definition}

Note that CSP polytopes can be defined independently of standard pairs
by the conditions of Corollary~\ref{cor:CSPchar}.

\begin{theorem}\label{thm:obstruction}  
  The following are equivalent. 
\begin{enumerate} 
\item The ideals $\ia$ and $\ica$ are not equal.
\item There is a generic $\omega \in \R^n$ for which $\A$ has a CSP
 polytope.
\item For every generic $\omega \in \R^n$, $\A$ has a CSP
 polytope.
\end{enumerate}
\end{theorem} 

\begin{proof} 
  The ideal $\ia = \ica$ if and only if for any generic $\omega \in
  \R^n$, $\inomega(\ia) = \inomega(\ica)$ which is if and only if
  $\inomega(\ica)$ has no CSPs.
\end{proof}

\section{Associated Primes of the Circuit Ideal}
\label{sec:assprimes}
In this section, we show how the associated primes of $\ica$ relate to
the CSP polytopes of its initial ideals discussed in
Section~\ref{sec:initialideals}. Recall that a face $F$ of $\cone(\A)$
is recorded as the set $\sigma := \{j \,:\, \a_j \in F\} \subseteq
[n]$.

\begin{proposition} \cite{EiSt} 
\label{prop:primarydecomposition}
\begin{enumerate}
\item The toric ideal $\ia$ is the radical of the circuit ideal $\ica$
  and hence the unique minimal prime of $\ica$. Furthermore, $ \ia
  = ( \ica \, : \, (x_1x_2 \cdots x_n)^{\infty})$.
\item All associated primes of $\ica$ are of the form $P_{\sigma} +
  \ia$ for some face $\sigma$ of $\cone(\A)$.  In particular, $\ia =
  P_{[n]} + \ia$. However, not all faces of $\cone(\A)$ need index an
  associated prime of $\ica$.
\end{enumerate}
\end{proposition}

\begin{remark}
\label{rem:otherideals}
  If $I = (\ica \, : \, \x^{\m})$ for some monomial $\x^{\m}$, then 
  $\ica \subseteq I \subseteq \ia$ and 
  Proposition~\ref{prop:primarydecomposition} holds for $I$.  
  Thus most of the results in this section stated for circuit ideals
  actually hold for all such ideals $I$.  
\end{remark}

\begin{definition} \cite{AtMc} 
\label{def:witness}
  Let $I$ be any ideal in $k[\x]$ and let $P$ be an ideal that
  contains $I$. Then $P$ is an associated prime of $I$ if 
  $P$ is prime and there exists some $f \in k[\x]$ such that $(I:f) \,
  = \, P$.  We call $f$ a \emph{witness} for $P$.
\end{definition}

Using Proposition~\ref{prop:primarydecomposition}, we can now state the 
main results of this section.  We say that $\tau$ is the {\em type} of a 
standard pair of the form $(\cdot, \tau)$. 

\begin{theorem}
\label{thm:primestoCSPs}
Let $\tau$ be any (possibly empty) proper face of $\cone(\A)$ and
$\omega$ be a generic weight vector.  If $P_{\tau} + \ia$ is
associated to $\ica$ (and hence embedded), then there exists a
circuit-specific standard pair of $\inomega(\ica)$ of type $\beta$ 
such that
  \begin{enumerate}
  \item $\beta \subseteq \tau$,
  \item if $\sigma$ is a face of $\cone(\A)$ properly contained in
    $\tau$, then $\beta$ is not contained in $\sigma$, and 
  \item $\left|{\beta}\right| = \dim \cone(\A_{\tau})$.
  \end{enumerate}
  Furthermore, there is a witness for the prime $P_{\tau} + \ia$ 
  whose leading term with respect to $\omega$ lies on such a CSP.
\end{theorem}

We also prove a partial converse. 

\begin{theorem} 
\label{thm:CSPstoprimes}
For a generic $\omega$, if $\inomega(\ica)$ has a CSP of type $\beta$,
then $\ica$ has an embedded prime $P_{\sigma} + \ia$ for some face
$\sigma$ of $\cone(\A)$ such that $\sigma \supseteq \beta$.
\end{theorem}

Before proving the theorems, we consider a few implications. We say
that a face $\tau$ of $\cone(\A)$ is {\em simplicial} if $| \tau | =
\dim \cone(\A_{\tau})$. If $\tau$ is a simplicial face of $\cone(\A)$,
then no binomial in $\ia$ is supported entirely on $\tau$, so
$P_{\tau} + \ia$ is just the monomial prime $P_{\tau} = \langle x_i: i
\notin \tau \rangle$.  Then Theorem~\ref{thm:primestoCSPs} specializes
as follows. 

\begin{corollary}
\label{cor:simplicialcase}   
If $P_{\tau} + \ia$ is an embedded prime of $\ica$ and $\tau$ is a
simplicial face of $\cone(\A)$, then for every generic $\omega$,
$\inomega(\ica)$ has a CSP of type $\tau$.
\end{corollary}

The situation is more complicated when non-simplicial faces of
$\cone(\A)$ index embedded primes.  In particular,
Theorem~\ref{thm:primestoCSPs} does not specify a particular $\beta
\subseteq [n]$ such that every monomial initial ideal of $\ica$ must
have a CSP of type $\beta$.

\begin{figure}
\begin{center}
\epsfig{file=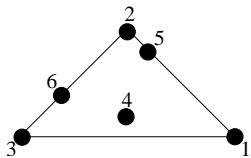,height=2.0cm} 
\end{center}
\caption{The point configuration of Example~\ref{ex:triangle}.}
\label{fig:triangle}
\end{figure}

\begin{example} 
\label{ex:triangle}
  Consider the matrix 
  $$A = \left( \begin{array}{rrrrrr}
    5  &  0 &  0 &  2 & 1 & 0 \\
    0  &  5 &  0 &  1 & 4 & 2 \\
    0  &  0 &  5 &  2 & 0 & 3 \end{array} \right).$$

The configuration $\A$ labeled $1, \ldots, 6$ in
Figure~\ref{fig:triangle} spans the cone over a triangle in $\R^3$, so
by Proposition~\ref{prop:primarydecomposition}, there are seven
possible embedded primes corresponding to the seven proper faces of
$\cone(\A)$.  All seven of these primes are indeed associated to $\ica$.  \
The two non-simplicial $2$-dimensional faces, $\{1,2,5\}$ and $\{2,3,6\}$,
index the primes $P_{\{1,2,5\}} + \ia = \langle f, d, c, ab^4-e^5
\rangle$ and $P_{\{2,3,6\}} + \ia = \langle~e,d,a,b^2c^3-f^5\rangle$.
The third $2$-face $\{1,3\}$ is simplicial and indexes the prime
$P_{\{1,3\}} = \langle b, d, e, f \rangle$.  The remaining four primes
$P_{\{1,2\}} = \langle c,d,e,f \rangle$, $P_{\{1,3\}} = \langle
b,d,e,f \rangle$, $P_{\{2,3\}} = \langle a,d,e,f \rangle$, and
$P_{\{\emptyset\}} = \langle a,b,c,d,e,f \rangle$ correspond to the
three rays of $\cone(\A)$ and to the apex, all of which are trivially
simplicial.

By Corollary~\ref{cor:simplicialcase}, each $\inomega(\ica)$ must have
CSPs of types $\{1,3\}$, $\{3\}$, $\{2\}$, $\{1\}$, and $\emptyset$
corresponding to the five simplicial faces of $\cone(\A)$. Since
$P_{\{a,b,e\}} + \ia$ is associated, Theorem~\ref{thm:primestoCSPs}
requires that $\inomega(\ica)$ has a CSP of type $\{1,2\}$, $\{1,5\}$,
or $\{2,5\}$.  Similarly, because of $P_{\{b,c,f\}} + \ia$,
there must always be a CSP of type $\{2,3\}$, $\{2,6\}$, or $\{3,6\}$.
We list the types of CSPs that appear for two term orders.

\begin{itemize} 
\item For lexicographic order with $f \succ e \succ \ldots \succ a$, 
  $\tin_{\succ}(\ica)$ has the following types of CSPs:
   $\{1,3\}, \{3\}, \{2\}, \{1\}, \emptyset, \{1,2\}, \{2,3\}$. 
\item For reverse lexicographic order with $a \succ b \succ \ldots
  \succ f$, $\tin_{\succ}(\ica)$ has the following types of CSPs:
$\{1,3\}, \{3\}, \{2\}, \{1\}, \emptyset, \{6\}, \{1, 5\}, \{3, 6\},
    \{2, 6\}, \{2, 5\} $.
\end{itemize}
\end{example}

We now prove Theorem~\ref{thm:primestoCSPs}.  The idea is
to find a witness for the embedded prime $P_{\tau} + \ia$, compute its
normal form with respect to the reduced Gr\"obner basis $G_\omega(\ica)$, 
and show that the result is a witness whose $\omega$-initial term lies 
on a CSP satisfying all of the desired properties.

\begin{lemma} \label{lem:witnesses}
  If $f$ is a witness for an embedded prime 
  $P_{\tau} + \ia$ of $\ica$, then the following hold.
\begin{enumerate} 
  \item The witness $f$ is in the toric ideal $\ia$.
  \item For any $g \in \ica$, $f+g$ is also a witness for $P_{\tau} +
    \ia$.  In particular, the normal form of $f$ with respect to
    $G_{\omega}(\ica)$ is a witness.
  \item If $\x^\m$ is a monomial with $\supp(\m) \subseteq \tau$, then 
    $\x^\m f$ is also a witness. 
\end{enumerate}
\end{lemma}

\begin{proof} 
\begin{enumerate}
\item Since $\tau$ is a proper face of $\cone(\A)$, there is some
  variable $x_i \in P_{\tau} + \ia$, so $x_i f \in \ica \subset \ia$.
  Since $\ia$ is a prime ideal without monomials, $f \in \ia$.
  \item Since $g \in \ica$, so is $pg$ for any polynomial $p \in \k[\x]$, 
    and thus $p(f+g) \in \ica \Leftrightarrow pf \in \ica$.  Thus 
    $(\ica \, : \, f+g) = (\ica \, : \, f) = P_{\tau} + \ia$. 
  \item   If $h \in P_{\tau} + \ia$, then $h (\x^\m f) = (\x^\m h)f$ 
    is in $\ica$ by the assumption that $f$ is a witness.  On the other 
    hand, if $h \notin P_{\tau} + \ia$, then neither is $\x^\m h$ because 
    $\supp(\m) \subseteq \tau$ and $P_{\tau} + \ia$ is prime.  Thus 
    $\x^\m h \notin P_{\tau} + \ia = (\ica \, : \, f)$,
    so $h \notin (\ica \, : \, \x^\m f)$.  Thus 
    $(\ica \, : \, \x^\m f) = (\ica \, : \, f)  = P_{\tau} + \ia$ as claimed. 
\end{enumerate}
\end{proof}

\begin{lemma}  \label{lem:firstwitness}
  If $f$ is a witness for an embedded prime $P_{\tau} + \ia$ of $\ica$
  and $\bar{f}$ is the normal form of $f$ with respect to
  $G_{\omega}(\ica)$, then $\inomega(\bar{f})$ lies on a CSP $(\cdot,
  \beta)$ of $\inomega(\ica)$ with $\beta \subseteq \tau$.
\end{lemma}

\begin{proof}
  By Lemma~\ref{lem:witnesses}, $\bar{f}$ is also a witness for
  $P_{\tau} + \ia$ and $\bar{f} \in \ia$.  This implies that
  $\inomega(\bar{f}) \in \inomega(\ia) \setminus \inomega(\ica)$, so
  $\inomega(\bar{f})$ must lie on some CSP $(\cdot, \beta) $of 
  $\inomega(\ica)$.  Since $x_i \bar{f} \in \ica$ whenever $i \notin \tau$
  because $\bar{f}$ witnesses $P_{\tau} + \ia$, it follows that $x_i
  \inomega(\bar{f}) \in \inomega(\ica)$ for $i \notin \tau$, so $\beta
  \subseteq \tau$.
\end{proof}

\vspace{.2cm} \smallskip\noindent {\it Proof of
  Theorem~\ref{thm:primestoCSPs}}: Suppose $P_{\tau} + \ia$ is an
embedded prime of $\ica$ and $e := \dim \cone(\A_{\tau})$.  We first
claim the following: there is a constant $C$ such that for all 
sufficiently large $N$, $P_{\tau} + \ia$ has at least $N^e$ witnesses whose
normal forms with respect to $G_{\omega}(\ica)$ have distinct leading 
terms, and each such leading term $\x^{\p}$ has the property that every
component of $\p$ is bounded above by $CN$.

Suppose the claim is true.  By Lemma~\ref{lem:firstwitness}, each such
monomial $\x^{\p}$ must lie on a CSP of type $\beta$ with $\beta
\subseteq \tau$.  Each such standard pair contains at most
$C^{\left|\beta\right|}(N+1)^{\left|\beta\right|}$ monomials $\x^\p$ such
that $p_i$ is bounded above by $CN$. Since there are only
finitely many standard pairs for $\inomega(\ica)$, all the standard
pairs of type $\beta$ with $|\beta| < e$ together cover only at most
$\bigo(N^{e-1})$ of the monomials which is not enough to contain the $N^e$
leading terms $\x^{\p}$.  Thus some of these leading terms must lie on
standard pairs $(\cdot, \beta)$ with $|\beta| \geq e$.  Since by
Corollary~\ref{cor:initialassprimes}, each $\beta$ is a face of the
triangulation $\Delta_{\omega}$ of $\A$, this is only possible if
$|\beta| = e$ and $\beta$ is not contained in any face $\sigma$ of
$\cone(\A)$ whose dimension is less than $e$.  These are exactly the
types of standard pairs in the conditions of
Theorem~\ref{thm:primestoCSPs}.

Now we prove the claim.  Since $P_{\tau}+\ia$ and $\ica$ are both
$\A$-homogeneous, there exists an $\A$-homogeneous witness $f$ for
$P_{\tau}+\ia$.  Set $\x^\u := \inomega(f)$.
Since $e = \dim \cone(\A_{\tau})$, we can find an $e$-subset $\alpha$ of
$\tau$ such that the columns of $A$ indexed by $\alpha$ are linearly
independent.  Thus if $\m_1 \neq \m_2$ are supported only on
$\alpha$, then $A \m_1 \neq A \m_2$.

Consider all polynomials of the form $\x^{\m}f$ where $0 \leq m_i < N$
for $i \in \alpha$ and $m_i = 0$ for $i \notin \alpha$.  Such a
polynomial is $\A$-homogeneous of $\A$-degree $A \m+A \u$ and so is its
normal form with respect to $G_{\omega}(\ica)$ since $\ica$ is an
$\A$-homogeneous ideal.  Thus the normal forms of these $N^e$
polynomials are all $\A$-homogeneous of different degrees, so in particular
they all have distinct leading terms.  Furthermore, by parts (2) and (3) of
Lemma~\ref{lem:witnesses}, each such normal form is a witness for
$P_{\tau}+\ia$.

It remains to establish that if $\x^{\p}$ is the leading term of one
of the normal forms, then each component of $\p$ is bounded by a
constant multiple of $N$.  Let $\a$ be a strictly positive vector in
the rowspan of $A$. Such a vector exists since $\la \cap \N^n =
\{0\}$. By scaling, we can assume that the minimum component of $\a$
is 1.  Let $R$ be its maximum component.  Since $A \p = A(\u + \m)$,
it follows that $\a \cdot \p = \a \cdot (\u + \m)$.  Then $\left
  \|\p\right \|_1 = \sum_{i=1}^n p_i \leq \sum_{i=1}^n a_i p_i = $
$$
\sum_{i=1}^n a_i(u_i + m_i) \\
\leq R \sum_{i=1}^n (u_i+m_i) = R(\sum_{i=1}^n u_i + \sum_{i=1}^n m_i)
< R(\|\u\|_1 + nN). $$
It follows that for any $i$, we have
  $$ p_i \leq \|\p\|_1 < RnN + R\|\u\|_1 $$ which is a bound of the desired
  form.
\endproof

We now prove Theorem~\ref{thm:CSPstoprimes}.  Recall the following
algebraic fact.

\begin{lemma} \label{lem:primesofquotient}
  If $I$ is an ideal in $k[\x]$ and $g$ is any polynomial, then the
  associated primes of $(I \, : \, g^{\infty})$ are exactly the
  associated primes of $I$ that do not contain $g$.
\end{lemma}

\begin{proposition} \label{prop:primesofsaturation}
  Recall that $\x_{\tau} = \prod_{i \in \tau} x_i$.  The
  associated primes of $(\ica \, : \, \x_{\tau}^{\infty})$ are exactly
  the primes $P_\sigma + \ia$ of $\ica$ that satisfy $\sigma \supseteq
  \tau$.
\end{proposition}

\begin{proof}
  We get $\x_{\tau} \in P_{\sigma} + \ia$ if and only if some variable $x_i$
  with $i \in \tau$ lies in $P_{\sigma} + \ia$, which in turn occurs if and
  only if $\tau$ is not contained in $\sigma$.  Now apply
  Lemma~\ref{lem:primesofquotient}.
\end{proof}
   
\medskip
\noindent {\it Proof of Theorem~\ref{thm:CSPstoprimes}}:
Suppose $(\x^{\u}, \beta)$ is a CSP of $\inomega(\ica)$.  Choose $f
\in \ia$ such that $\inomega(f) = \x^{\u} \x_{\beta}^{\m}$ for some
$\m \geq 0$. This is possible since every CSP of $\inomega(\ica)$
contains non-standard monomials of $\inomega(\ia)$. Since no monomial
of the form $\inomega(f) \cdot \x_{\beta}^\ast$ lies in
$\inomega(\ica)$, no polynomial of the form $f \cdot \x_{\beta}^\ast$
lies in $\ica$. This implies that $(\ica \, : \, \x_{\beta}^{\infty})$
does not contain $f$ and is hence not equal to $\ia$.  However, since
$(\ica \, : \, x_{\beta}^{\infty}) \subsetneq \ia$, 
$(\ica \, : \, x_{\beta}^{\infty})$ must have an embedded prime.  This 
prime is also embedded in $\ica$ by Proposition~\ref{prop:primesofsaturation}, 
and it has the form $P_\sigma + \ia$ for some $\sigma \supseteq \beta$. 
\endproof

Theorem~\ref{thm:CSPstoprimes} is only a partial converse to
Theorem~\ref{thm:primestoCSPs}.  It is not true for a given weight
vector $\omega$ that the existence of a CSP $(\x^{\u}, \beta)$ of
$\inomega(\ica)$ satisfying the conditions of
Theorem~\ref{thm:primestoCSPs} with respect to some proper face $\tau$
of $\cone(\A)$ implies that $P_{\tau} + \ia$ is associated to $\ica$.

\begin{example} \label{ex:falseempty}
Take
$$A = \left( \begin{array}{rrrr}
    1 & 3 & 2 & 4 \\
    1 & 4 & 5 & 2 \end{array} \right).$$
The values of the $\A$-graded Hilbert
function of this $\A$ are shown in Figure~\ref{fig:hilbvalues}.  The
proper faces of $\cone(\A)$ are $\{3\}$, $\{4\}$, and $\emptyset$.
Only the first two index associated primes of $\ica$.  However,
if we take $\omega$ to represent the lexicographic term order with $a
\succ b \succ c \succ d$, there are five CSPs of $\inomega(\ica)$ of
type $\emptyset$.  On the other hand, if $\omega$ represents the $\A$-graded
reverse lexicographic order with $a \succ b \succ c \succ d$, then
there are no CSPs of type $\emptyset$.
\end{example}

Based on Example~\ref{ex:falseempty}, it is possible that if $\tau$ is 
a face of $\cone(\A)$ such that for {\bf every} generic $\omega$ there 
is a CSP of $\inomega(\ica)$ satisfying the conditions of 
Theorem~\ref{thm:primestoCSPs} with respect to $\tau$, then $P_{\tau} + \ia$ 
is associated to $\ica$.

We conclude this section with an application of
Theorem~\ref{thm:primestoCSPs}. Let $\A$ be a configuration satisfying
$\Z \A = \Z^d$ and whose vectors comprise the lattice points in a
lattice polytope $R$.  Further assume that the first row of $A$ is
$(1,\ldots,1)$, so $R$ is at height 1. The polytope $R$ defines a
projective toric variety $\xa$ and the faces $\{\tau\}$ of $R$ (which
are in bijection with the faces of $\cone(\A)$) index a canonical
collection of affine charts $\{U_{\tau}\}$ covering $\xa$ \cite{Ful}.
We investigate how smoothness of $U_{\tau}$ determines whether
$P_{\tau} + \ia \in \textup{Ass}(\ica)$.

\begin{definition} 
\begin{enumerate} 
\item Let $K$ be a convex rational polyhedral cone in $\R^t$ that does
  not contain a line.  We say that $K$ is \emph{smooth} if it is
  generated by primitive vectors that form part of a basis for $\Z^t$.
\item Let $K_F$ denote the inner normal cone of a face $F$ of a
  polytope $Q$. The face $F$ is \emph{smooth} if the restriction 
  of $K_F$ to the linear span of $Q$ is smooth.
\end{enumerate}
\end{definition}  

\begin{remark} 
\begin{enumerate}
\item If $v$ is a smooth vertex of a polytope $Q$ then there are
  exactly $\dim \, Q$ edges of $Q$ incident to $v$. Further, the cone
  dual to $K_v$ is also smooth \cite[Theorem 2.10, Chapter V]{Ewald}.
  Note that this dual cone is the tangent cone of $Q$ at $v$ and
  contains $Q$.

\item A face $F$ of a polytope $Q$ is smooth if and only if the affine 
toric variety $U_F$ is smooth \cite{Ful}. 

\end{enumerate}
\end{remark}

\begin{theorem}\label{thm:smoothvertex} 
  Let $\A$ and $R$ be as above.  If $\a_n$ is a smooth vertex of $R$, then 
  $P_{\{n\}} (= P_{\{n\}} + \ia)$ is not an associated prime of $\ica$.
\end{theorem}

\begin{proof}   Suppose $P_{\{n\}}$ is associated. Since $\{n\}$ is a
  simplicial face of $\cone(\A)$, by
  Corollary~\ref{cor:simplicialcase}, every monomial initial ideal 
  $\inomega(\ica)$ has a CSP of the form $(\x^{\u}, \{n\})$. In 
  particular, let $\omega$ represent an elimination order with $x_n$ most 
  expensive.  We may assume that each of $\a_1, \ldots, \a_{d-1} \in \A$ 
  is the first lattice point from $\a_n$ along one of the $d-1$ edges 
  incident to $\a_n$. Let $\y_i := \a_i - \a_n$ for $i = 1, \ldots, d-1$. 
  
  Since $\a_n$ is smooth and $R$ is contained in the tangent cone at $\a_n$, 
  for each lattice point in $R$ (i.e. column of $A$), there are unique 
  $m_i \in \N$ such that $\a_j = \a_n
  + \sum_{i=1}^{d-1} m_i \y_i.$ Rearranging terms, and setting $M = -1
  + \sum_{i=1}^{d-1} m_i$, we get $\a_j + M \a_n = \sum_{i=1}^{d-1}
  m_i \a_i$ with all coefficients nonnegative.  If $j = n$, this
  equation reduces to $0 = 0$, and if $1 \leq j \leq d-1$, it reduces to
  $\a_j = \a_j$.  But in the nontrivial case where $d-1 < j < n$, this
  relation is a \emph{circuit} because $a_1, \ldots, a_{d-1}, a_n$ form 
  a maximal linearly independent set.
  Thus $x_j x_n^M - \prod_{i=1}^{d-1} x_i^{m_i} \in \ica$.  By choice of
  $\omega$, its leading term is $x_j x_n^M$.  Since $(\x^{\u}, \{n\})$ 
  is a CSP, this term must not divide $x_n^N \x^{\u}$ for any $N$.  This 
  means that $j$ is not in the support of $\u$.
  
  For $N$ sufficiently large, $x_n^N \x^{\u} \in \inomega(\ia)$, so we
  can choose $\x^\v \not \in \inomega(\ia)$ such that $x_n^N \x^{\u} -
  \x^{\v} \in \ia$.  Since $\ia$ is prime, factor out any common
  monomial to get $\x^{\tilde{\u}} - \x^{\tilde{\v}} \in \ia$ where
  $\tilde{\u}$ and $\tilde{\v}$ have disjoint supports.  Since
  $\tilde{\u} - \tilde{\v} \in \la$, the convex hulls of $\{\a_i \,:\,
  i \in \supp(\tilde{\u}) \}$ and $\{\a_i \,:\, i \in
  \supp(\tilde{\v}) \}$ must intersect.
  
  Since $\a_n$ is smooth, we can assume by applying an invertible
  $\Z$-affine transformation that $\a_n$ is the origin and $\a_i$ is
  the $i$th standard basis vector in $\Z^{d-1}$ for $1 \leq i
  \leq d-1$.  That is, $\a_n, \a_1, \a_2, \ldots, \a_{d-1}$ are the
  vertices of the standard simplex $S$ in $\R^{d-1}$.
 
  Since $j \notin \supp(\tilde{\u})$ for any $d-1 < j < n$, $U :=
  \conv(\a_i \,:\, i \in \supp(\tilde{\u}))$ is a face of $S$.  Since
  $\supp(\tilde{\v}) \cap \supp(\tilde{\u}) = \emptyset$ and $S$
  contains no lattice points except its vertices, $V := \conv( \a_i
  \,:\, i \in \supp(\tilde{\v}) )$ consists only of vertices of $S$
  outside $U$ along with lattice points in $R \setminus S$.  Now $S$
  and $\overline{R \setminus S}$ are both convex, so $U$ and $V$ could
  intersect only on the boundary of $S$.  But since the vertices in
  $U$ and those in $V \cap S$ form disjoint faces of $S$, there is no
  intersection on this boundary, contradicting that $U \cap V \neq
  \emptyset$.
\end{proof}

\begin{example} 
  {\bf Non-smooth vertices of $R$ may or may not index associated
    primes of $\ica$}: For the following $A$, the polytope $R$
  is a triangle in $\R^3$.
  $$
  A = \left( \begin{array}{rrrrr}
      1 & 1 & 1 & 1 & 1 \\
      0 & 3 & 4 & 5 & 5 \\
      0 & 1 & 1 & 1 & 2 \end{array} \right).$$
  None of the three vertices $(0, 0)$, $(5, 1)$, and $(5, 2)$ of $R$ are
  smooth.  The vertices $(0, 0)$ and $(5,2)$ both index associated primes, 
  while the vertex $(5,1)$ does not.
\end{example}

\begin{example} {\bf Smooth edges of $R$ may index associated primes of
    $\ica$}: Consider the matrix
  $$A = \left( \begin{array}{rrrrrrrrrrrrrrrrrr}
    1 & 1 & 1 & 1 & 1 & 1 & 1 & 1 & 1 & 1 & 1 & 1 & 1 & 1 & 1 & 1 & 1 & 1 \\
    0 & 0 & 0 & 0 & 0 & 0 & 0 & 0 & 0 & 1 & 1 & 1 & 1 & 1 & 1 & 1 & 1 & 1 \\
    0 & 0 & 0 & 1 & 1 & 1 & 2 & 2 & 2 & 0 & 0 & 0 & 1 & 1 & 1 & 2 & 2 & 2 \\
    0 & 1 & 2 & 0 & 1 & 2 & 0 & 1 & 2 & 0 & 1 & 2 & 0 & 1 & 2 & 0 & 1 & 2 
\end{array} \right)$$  
where $\cone(\A)$ is the cone over a rectangular prism $R$.  All edges
of $R$ are smooth, but we compute that four of them index
associated primes of $\ica$.
\end{example}

\section{Fans of Toric and Circuit Ideals}
\label{sec:fans}
The main goal of this section is to compare $\ia$ and $\ica$ via three
polyhedral fans that can be associated to them. These results rely on
Theorem~\ref{thm:fatstdpairsagree} which states that for a generic
$\omega$, the fat standard pairs of $\inomega(\ia)$ and
$\inomega(\ica)$ are the same. We begin by providing a second proof of
Theorem~\ref{thm:fatstdpairsagree} using completely different
techniques from those used in Section~\ref{sec:initialideals}. These
methods exploit the general structure of standard pairs and can be
applied to any pair of ideals, one contained in the other, that
satisfy hypotheses similar to $\ia$ and $\ica$.  Recall the notions of
\emph{admissible pair} and \emph{standard pair} from
Section~\ref{sec:initialideals}.

\begin{definition} \label{pair}
  A \emph{pair} $(\x^\u,\sigma)$ is the set of monomials $\{\x^\m:
  \m\geq\u \textup{ and } m_i=u_i \textup{ for all } i\not\in\sigma\}
  \subset \k[\x]$. The \emph{dimension} of the pair $(\x^\u,\sigma)$ is
  the cardinality of $\sigma$.
\end{definition}

\begin{lemma}\cite{CLO}
\label{lem:ddimdpair}
  Let $M\subseteq \k[\x]$ be a monomial ideal with the Krull dimension of 
  $\k[\x]/M$ equal to $d$.  Then any admissible pair of $M$ has dimension 
  at most $d$.
\end{lemma}

\begin{lemma} \label{lem:dimcircuitpair}
  Let $M,M'\subseteq\k[\x]$ be monomial ideals with $M\subseteq M'$
  and with the Krull dimension of $\k[\x]/M$ equal to $d$. Then any
  $d$-dimensional standard pair $(\x^\u,\sigma)$ of $M'$ is also a
  standard pair of $M$.  In particular, any $d$-dimensional (fat)
  standard pair for $\inomega(\ia)$ is also a standard pair for
  $\inomega(\ica)$.
\end{lemma}

\begin{proof}
  The admissible pair $(\x^\u,\sigma)$ for $M'$ is admissible for $M$
  since $M\subseteq M'$. To prove that the pair is standard for $M$ we
  must show that there does not exist another admissible pair
  $(\x^\v,\tau)$ of $M$ with $\x^\v|\x^\u$ and
  $\supp(\u-\v)\cup\sigma\subseteq\tau$.  Suppose such a $(\x^\v,
  \tau)$ exists.  Since $\dim(\k[\x]/M)=d$ and $M\subseteq M'$,
  $\dim(\k[\x]/M') \leq d$, so by Lemma~\ref{lem:ddimdpair},
  $|\tau|\leq d$.  Then $\sigma\subseteq \tau$, $|\sigma|=d$, and
  $|\tau|\leq d$ implies $\sigma=\tau$. Since
  $\supp(\u)\cap\sigma=\emptyset$, $\supp(\u-\v) \cap \sigma$ must
  also be empty. On the other hand we have $\supp(\u-\v)\subseteq
  \tau=\sigma$. Hence $\supp(\u-\v)=\emptyset$, so $\u=\v$.
\end{proof}

\begin{definition}
  Fix $\a\in\N_{>0}^n$ and grade $\k[\x]$ by $\a$ via $\deg(x_i) :=
  a_i$.  Let $I \subseteq \k[\x]$ be an $\a$-homogeneous ideal and let
  $I_s$ be the $\k$-vector space of polynomials in $I$ of $\a$-degree
  $s$. Let $H_{I,\a}(s) := \dim_\k({\k[\x]_s / I_s})$ be the
  $\a$-\emph{graded Hilbert function} of $I$.
\end{definition}

If $I$ is $\a$-homogeneous, then $H_{I,\a}(s) = H_{\inomega(I),\a}(s)$
for any $\omega \in \R^n$.  If $I$ is a monomial ideal, $H_{I,\a}(s)$
equals the number of standard monomials of $I$ of degree $s$.

Let $(\x^\u,\sigma)_s$ denote the set of monomials of $\a$-degree $s$
on the pair $(\x^\u,\sigma)$.  Hence
$(\x^\u,\sigma)=\bigcup_{s=0}^\infty (\x^\u,\sigma)_s$.  We use the
notation $\bigo (s^e)$ for the set of all functions $f : \N
\rightarrow \Z$ such that $|f(s)|$ is bounded above by a polynomial in
$s$ of degree $e$.  Note that $\bigo (s^e)$ is closed under addition.

\begin{lemma}
\label{lem:subspacegrowth}
Let $(\x^\u,\sigma)$ be a pair and define the function
$f_\a(s)=|(\x^\u,\sigma)_s|$. Then
$f_\a\in\bigo(s^{|\sigma|-1})\backslash\bigo(s^{|\sigma|-2})$.
\end{lemma}

\begin{proof}
  We will show this for $\u=0$.  The generalization of the result is
  straightforward. Let ${\bf 1} := (1, \dots,1)^T\in\N^n$ and $b$ be
  the product of the entries in $\a$.
  
  If $\x^\v\in(1,\sigma)$ with $\a$-degree $s$, then $\x^{(v_1
    a_1,\dots,v_n a_n)}$ has total degree $s$ and is also in
  $(1,\sigma)$.  Hence $f_\a(s)\leq f_{\bf 1}(s)$.  On the other hand,
  if $\x^\v\in(1,\sigma)$ with total degree $s$, then
  $\x^{({\frac{b}{a_1}}v_1,\dots,{\frac{b}{a_n}} v_n)}$ is in
  $(1,\sigma)$ and has $\a$-degree $bs$. Hence $f_{\bf 1}(s)\leq
  f_\a(bs)$.
  
  Now$f_{\bf 1}(s) = \binom{s+|\sigma|-1}{|\sigma|-1}$ is a polynomial
  of degree $|\sigma|-1$. The fact that $f_\a(s)\leq f_{\bf 1}(s)$
  implies that $f_\a\in\bigo(s^{|\sigma|-1})$. If
  $f_\a\in\bigo(s^{|\sigma|-2})$, then so is the function $s\mapsto
  f_\a(bs)$.  Then by the previous paragraph, so is $f_{\bf 1}$ which
  is a contradiction.
\end{proof}

\begin{lemma}
\label{lem:intersectpairs}
Either two pairs $(\x^\u,\sigma)$ and $(\x^\v,\tau)$ do not intersect or 
$(\x^\u,\sigma)\cap(\x^\v,\tau) = (\x^{\max(\u,\v)},\sigma\cap\tau)$.
\end{lemma}

\begin{proof}
  Assume that $(\x^\u,\sigma)$ and $(\x^\v,\tau)$ have a non-trivial
  intersection.

  Let $\x^\m$ be any monomial in $(\x^\u,\sigma)\cap(\x^\v,\tau)$.
  Then $\m\geq\max(\u,\v)$ since $\m\geq\u$ and $\m\geq\v$. Let
  $i\in\overline{(\sigma\cap\tau)}
  =\overline{\sigma}\cup\overline{\tau}$. Without loss of generality
  $i\in\overline{\sigma}$ which implies that $m_i=u_i$ since $\x^\m
  \in (\x^\u, \sigma)$. Therefore, $u_i \geq v_i$ since
  $\m\geq\max(\u,\v)$. This implies that $m_i=\max(\u,\v)_i$, and
  $(\x^\u,\sigma)\cap(\x^\v,\tau) \subseteq
  (\x^{\max(\u,\v)},\sigma\cap\tau)$.
  
  Choose $\x^\m \in (\x^\u,\sigma)\cap(\x^\v,\tau)$ not divisible by
  any other monomial in the intersection.  Then we have $\x^\u
  \x^\alpha = \x^\m = \x^\v \x^\beta$ where $\supp(\alpha) \subseteq
  \sigma$ and $\supp(\beta) \subseteq \tau$. Further, $\supp(\alpha)
  \cap \supp(\beta) = \emptyset$ since $\x^\m$ is not divisible by any
  other monomial in $(\x^\u,\sigma)\cap(\x^\v,\tau)$.  Then for each
  $i \in [n]$, since $i$ cannot occur in both $\supp(\alpha)$ and
  $\supp(\beta)$, $m_i$ = $\max(u_i, v_i)$.  Thus $\m = \max(\u,\v)$,
  so $(\x^\m, \sigma\cap\tau) = (\x^{\max(\u,\v)},\sigma\cap\tau)
  \subseteq (\x^\u,\sigma)\cap(\x^\v,\tau).$
\end{proof}

\begin{figure}[ht]
\begin{center}
\epsfig{file=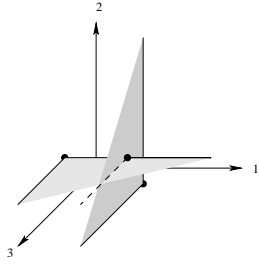,height=3.4cm} 
\end{center}
\caption{ The intersection of the pairs $(x_1^3x_3,\{2,3\})$ and 
$(x_2^2x_3^2,\{1,3\})$ is the pair $(x_1^3x_2^2x_3^2,\{3\})$.}
\end{figure}

\begin{lemma}
\label{lem:intersectstdpairs}
Let $M \subset \k[\x]$ be a monomial ideal. Let $(\x^\u,\sigma)$ and
$(\x^\v,\tau)$ be different standard pairs of $M$ with non-empty
intersection $(\x^{\max(\u,\v)},\sigma\cap\tau)$.  Then
$\sigma\cap\tau$ is properly contained in each of $\sigma$ and $\tau$.
\end{lemma}
\begin{proof}
  Suppose $\sigma\cap\tau=\sigma$.  Choose $\x^\m$ in the intersection
  of the two standard pairs.  If $i \in \sigma$, $v_i=u_i=0$ since
  $\sigma\subseteq\tau$ and $(\x^\u,\sigma)$ and $(\x^\v,\tau)$ are
  admissible. If $i\in\overline{\sigma}$, we have $m_i=u_i$ since
  $\x^\m\in(\x^\u,\sigma)$ and $u_i = m_i \geq v_i$ since
  $\x^\m\in(\x^\v,\tau)$.  In both cases $u_i\geq v_i$, implying
  $\u=\max(\u,\v)$. Consequently
  $(\x^\u,\sigma)=(\x^{\max(\u,\v)},\sigma)=(\x^{\max(\u,\v)},\sigma\cap\tau)
  =(\x^\u,\sigma)\cap(\x^\v,\tau)\subseteq(\x^\v,\tau)$. This implies
  that $(\x^\u,\sigma) = (\x^\v,\tau)$ since both are standard pairs,
  which is a contradiction.
\end{proof}

\begin{lemma}
\label{lem:dpairsandhilbertfunction}
Let $M, M'\subset \k[\x]$ be monomial ideals with $M\subseteq M'$
and let $d$ be the Krull dimension of $\k[\x]/M$.  Let $\a\in\N_{>0}^n$. 
Then the $d$-dimensional standard pairs of $M$ and $M'$ are the same if 
and only if $H_{M,\a}-H_{M',\a}\in \bigo(s^{d-2})$.
\end{lemma} 

\begin{proof}
  $(\Rightarrow)$: We count the number of standard monomials of each
  monomial ideal via the Inclusion-Exclusion Principle, getting a
  formula of the form 
  $$
  |(\x^{\alpha_1},\sigma_1)_s\cup \dots
  \cup(\x^{\alpha_t},\sigma_t)_s| = \sum_{r=1}^t(-1)^{r-1}\sum_{1\leq
    i_1 < \dots < i_r\leq t} |(\x^{\alpha_{i_1}},\sigma_{i_1})_s\cap
  \dots \cap(\x^{\alpha_{i_r}},\sigma_{i_r})_s| $$
  for each of $M$ and
  $M'$.  The nonzero terms having $r\not=1$ or having
  $|\sigma_{i_1}|\not=d$ belong to $\bigo(s^{d-2})$ by Lemma
  \ref{lem:intersectstdpairs} and Lemma \ref{lem:subspacegrowth}. The
  remaining ones correspond to the $d$-dimensional standard pairs.
  Since the $d$-dimensional standard pairs of $M$ and $M'$ are the
  same, these remaining terms are the same for both sums.  Hence the
  difference $H_{M,\a}-H_{M',\a}$ is in $\bigo(s^{d-2})$.
  
  $(\Leftarrow)$: By Lemma \ref{lem:dimcircuitpair}, each term in
  the sum for $M'$ that is not in $\bigo(s^{d-2})$ is also in the sum 
  for $M$.  Suppose there were more terms in the sum for $M$ not in
  $\bigo(s^{d-2})$.  These terms would all be positive, so the difference
  $H_{M,\a}-H_{M',\a}$ would not be in $\bigo(s^{d-2})$, a contradiction. 
\end{proof}

\begin{lemma}
\label{lem:onecolon}
Let $I\subset\k[\x]$ be an $\a$-homogeneous ideal such that $\ica
\subseteq I\subseteq \ia$.  Then
$H_{I,\a}-H_{(I:x_1^\infty),\a}\in\bigo(s^{d-2})$.
\end{lemma}

\begin{proof}
  Let $\omega$ represent the $\a$-graded reverse lexicographic term
  order with $x_1 \prec x_2 \prec \dots \prec x_n$.  Since
  $I\subseteq(I:x_1^\infty)$ implies $\inomega(I)\subseteq
  \inomega(I:x_1^\infty)$, it suffices (by
  Lemma~\ref{lem:dpairsandhilbertfunction}) to show equality between
  the sets of $d$-dimensional standard pairs of these initial ideals
  to prove the lemma. Since $\ica \subseteq I \subseteq \ia$, all
  ideals involved are $d$-dimensional.
  
  Let $(\x^\u,\sigma)$ be a standard pair of $\inomega(I)$ with
  $|\sigma|=d$. Then by Remark~\ref{rem:Jstdpairs}, $(1,\sigma)$ is
  also a standard pair of $\inomega(I)$. Recall that the
  $\sigma$-columns of $A$ are independent over $\Q$ since $\sigma$ is
  a maximal face of $\Delta_{\omega}$.  We claim that $1\in\sigma$.
  If not, then $a_1$ and $\{a_i\}_{i\in\sigma}$ are dependent over
  $\Q$. By Lemma~\ref{lem:circuitscaling} this gives an
  $\a$-homogeneous binomial $\x^\alpha-\x^\beta \in \ica \subseteq I$
  with $\supp(\alpha), \supp(\beta)\subseteq\sigma\cup\{1\}$ and
  $\supp(\alpha) \cap \supp(\beta) = \emptyset$. The initial term
  $\inomega(\x^\alpha-\x^\beta)$ must be the term not containing $x_1$
  by the definition of $\omega$, so the support of this term is
  contained in $\sigma$. This contradicts $(1,\sigma)$ being standard.

  The ideal $\inomega(I:x_1^\infty)$ can be computed from $\inomega(I)$ by
  projection along the $x_1$ axis, see \cite[Lemma 12.1]{GBCP}. Since the
  monomials in $(\x^\u,\sigma)$ are all standard for
  $\inomega(I)$ with $1\in\sigma$, the pair $(\x^\u,\sigma)$ must be 
  admissible for $\inomega(I:x_1^\infty)$.  It is standard because its root 
  is 1 and its dimension is $d$.
\end{proof}

\noindent{\em Second proof of Theorem~\ref{thm:fatstdpairsagree}}: Let $\a$ be 
a positive vector in the rowspan of $A$.  By
Proposition~\ref{prop:primarydecomposition}, we have $\ia = (\ica \, :
\, (x_1 x_2 \cdots x_n)^{\infty}) = ((\ica \, : \, x_1^{\infty}) \,
: \ldots \, : x_n^{\infty})$ so by repeatedly applying
Lemma~\ref{lem:onecolon}, we see that $H_{\ia,\a} - H_{\ica, \a} \in
\bigo(s^{d-2})$.  It then follows by
Lemma~\ref{lem:dpairsandhilbertfunction} that $\inomega(\ia)$ and
$\inomega(\ica)$ have the same $d$-dimensional standard pairs.
\endproof

\subsection{Polyhedral Fans of \ica}
An ideal in $\k[\x]$ gives rise to several natural equivalence
relations on $\R^n$ some of which give rise to polyhedral fans. In
this final part, we compare various equivalence relations and fans for
toric and circuit ideals.

\begin{definition} \cite[page 112]{SST}
  Let $J\subset \k[\x]$ be an ideal and $d$ be the Krull dimension of 
  $\k[\x]/J$.
  Define $\top(J)$ to be the intersection of all primary components of
  $J$ of dimension $d$.
\end{definition}

Note that $\top(J)$ is well defined since the $d$-dimensional primary
components of $J$ are not embedded and are hence unique. 

\begin{definition}
\label{def:equivalence}
Let $I\subset \k[\x]$ be an ideal homogeneous with respect to a 
positive vector $\a\in\N_{>0}^n$. We define three equivalence 
relations on $\R^n$:
\begin{itemize}
\item The \emph{initial ideal} equivalence relation  
$\u\sim \v \Leftrightarrow \inu(I)=\inv(I)$.
\item The \emph{top} equivalence relation  
$\u\sim \v \Leftrightarrow \top(\inu(I))=\top(\inv(I))$.
\item The \emph{radical} equivalence relation
$\u\sim \v \Leftrightarrow \sqrt{\inu(I)}=\sqrt{\inv(I)}$.
\end{itemize}
\end{definition}
In all three cases, the equivalence classes are invariant 
under translation by $\a$.

\begin{definition}
A collection $C$ of polyhedra in $\R^n$ is a \emph{polyhedral complex} if:
\begin{enumerate}
\item all proper faces of a polyhedron $P\in C$ are in $C$, and
\item the intersection of any two polyhedra $A,B\in C$ is a face of 
$A$ and $B$. 
\end{enumerate}
A polyhedral complex is a \emph{fan} if it only consists of cones.
\end{definition}
We say that an equivalence relation defines a fan $F$ if the closures 
of its equivalence classes are exactly the cones in $F$.

\begin{proposition}
Let $I$ be as in Definition \ref{def:equivalence}. Then
\begin{enumerate}
\item The initial ideal equivalence relation defines the 
\emph{Gr\"obner fan} of $I$.
\item The radical equivalence relation does not define a fan in general.
\end{enumerate}
\end{proposition}
A proof of the first claim is given in \cite[Chapter 2]{GBCP}. See
also \cite{MR}. The following example demonstrates the second claim.

\begin{figure}
\begin{center}
\epsfig{file=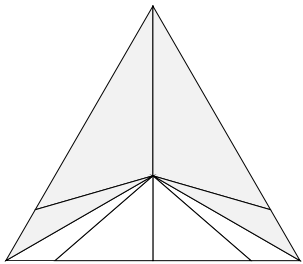,height=3.4cm}
\end{center}
\caption{The Gr\"obner fan from Example \ref{ex:radicalfan} and the
  two radical equivalence classes. The fan is drawn in the standard
  simplex with $(1,0,0)$ at the right bottom, $(0,1,0)$
  at the left bottom and $(0,0,1)$ at the top.}
\label{fig:radicalfan}
\end{figure}

\begin{example}
\label{ex:radicalfan}
The radical equivalence classes of the homogeneous ideal $I=\langle
c^4-ba^3,ab^3-ba^3\rangle\subset\k[a,b,c]$ are not all convex. This
ideal has eight monomial initial ideals. Four of them have radical
$\langle ab, ac, bc\rangle$ and the other four have radical $\langle
ab, c\rangle$. The intersection of the Gr\"obner fan with the
$2$-dimensional standard simplex is shown in Figure
\ref{fig:radicalfan} and the two radical equivalence classes are shown
in gray and white.
\end{example}

However, for toric ideals, all three equivalence relations of
Definition \ref{def:equivalence} give rise to polyhedral fans.

\begin{proposition}
\label{prop:toricrefinement}
\begin{enumerate}
\item The radical equivalence relation of $\ia$ defines the 
\emph{secondary fan} of $\A$.
\item The top equivalence relation of $\ia$ defines the 
\emph{hypergeometric fan} of $\A$.
\end{enumerate}
Furthermore, the Gr\"obner fan of $\ia$ is a refinement of the
hypergeometric fan of $\A$, which is a refinement of the secondary
fan of $\A$.

\end{proposition}
Proposition \ref{prop:toricrefinement} may be taken as the definition
of the hypergeometric and secondary fans of $\A$. The proposition is
a collection of several known results \cite[Proposition~3.3.1 and
Corollary~3.3.2]{SST}, \cite{BFS}, and \cite[Chapter~8]{GBCP}.
We now study the three equivalence classes for $\ica$.

\begin{theorem}
\label{thm:circuitrefinement}
The radical equivalence classes of $\ica$ form a polyhedral fan that
coincides with the secondary fan of $\ia$. 
\end{theorem}

\begin{proof}
  By Proposition \ref{prop:toricrefinement}, to prove that the radical
  equivalence relation of $\ica$ defines the secondary fan of $\A$ it
  suffices to show that $\sqrt{\inomega(\ia)}=\sqrt{\inomega(\ica)}$
  for any $\omega$. Since $\ica\subseteq \ia$ one inclusion is clear.
  To prove the other inclusion we first observe $\sqrt{\inomega(\ia)}$
  is $\omega$-homogeneous since $\inomega(\ia)$ is.  Hence it suffices
  to show that any homogeneous element in $\sqrt{\inomega(\ia)}$ is
  also in $\sqrt{\inomega(\ica)}$. Let $f\in \sqrt{\inomega(\ia)}$ be
  $\omega$-homogeneous. Then there exists some $m$ such that $f^m\in
  \inomega(\ia)$. The polynomial $f^m$ is also $\omega$-homogeneous,
  so $f^m=\inomega(F)$ for some $F\in\ia$.  Since $\sqrt{\ica} = \ia$,
  $F^{k}\in\ica$ for some $k$, and
  $\inomega(F^{k})=\inomega(F)^{k}=f^{mk}$.  Hence, $f\in
  \sqrt{\inomega(\ica)}$.
\end{proof}

Now we consider the top equivalence relation of $\ica$.

\begin{lemma}\cite[Lemma~3.2.4]{SST}
\label{lemma:topdependsonfat}
Let $M\subseteq\k[\x]$ be a monomial ideal such that the Krull 
dimension of $\k[\x]/M$ is $d$. Then
$$\top(M)=\bigcap_{(\x^\u,\sigma)}\langle x_i^{\u_i+1}:
i\not\in\sigma\rangle$$ where the intersection is taken over all 
$d$-dimensional standard pairs of $M$.
\end{lemma}

\begin{lemma}
\label{lem:circuittop}
For any generic $\omega\in\R^n$, 
$\top(\inomega(\ia))=\top(\inomega(\ica))$.
\end{lemma}
\begin{proof}
  For a generic $\omega$, the initial ideals $\inomega(\ia)$ and
  $\inomega(\ica)$ are both $d$-dimensional monomial ideals, and by
  Lemma \ref{lemma:topdependsonfat}, their $\top$s only depend on
  their $d$-dimensional standard pairs. By
  Theorem~\ref{thm:fatstdpairsagree}, these standard pairs are the
  same for both initial ideals.
\end{proof}

Lemma~\ref{lem:circuittop} proves the following.

\begin{proposition} \label{prop:top}
The top equivalence relation defines the same maximal cells for both
$\ia$ and $\ica$.  These cells are precisely the open maximal cones
in the hypergeometric fan of $\ia$.
\end{proposition} 

In contrast to Theorem~\ref{thm:circuitrefinement} and
Proposition~\ref{prop:top}, we have the following result for the
initial ideal equivalence relation for $\ia$ and $\ica$.

\begin{figure}[h]
\begin{center}
\epsfig{file=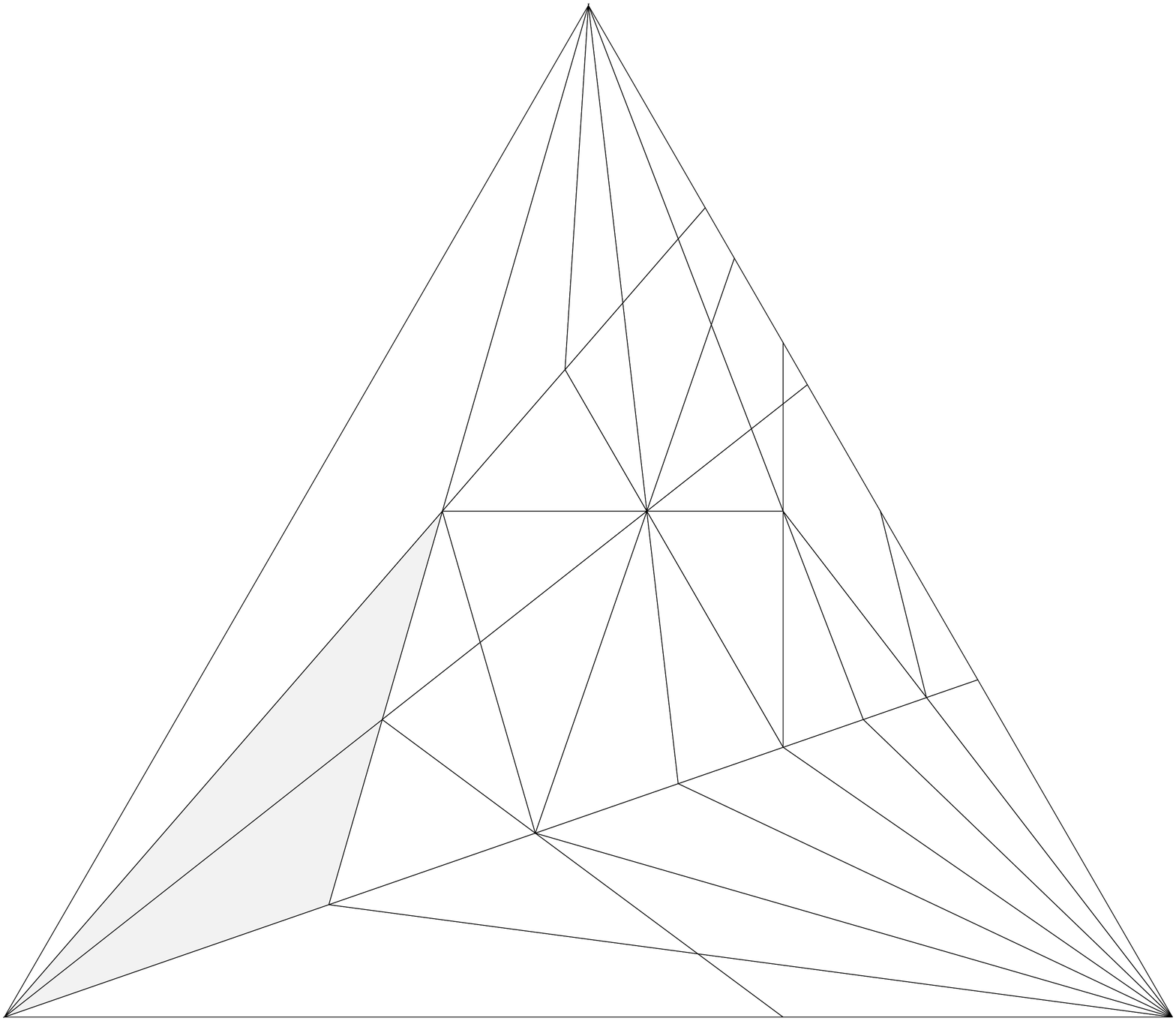,height=3.4cm} 
\hspace{0.01cm}
\epsfig{file=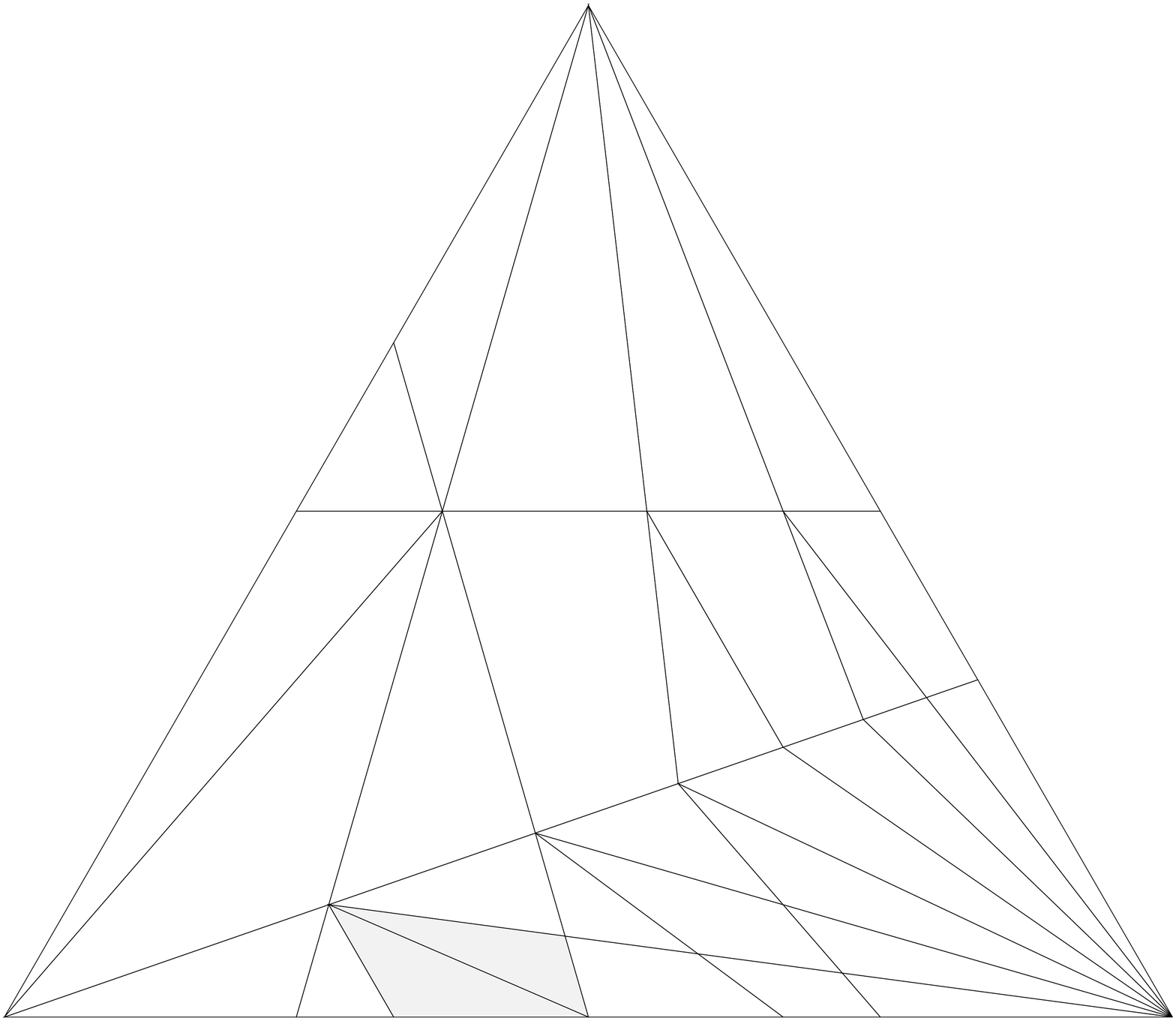,height=3.4cm} 
\hspace{0.01cm}
\epsfig{file=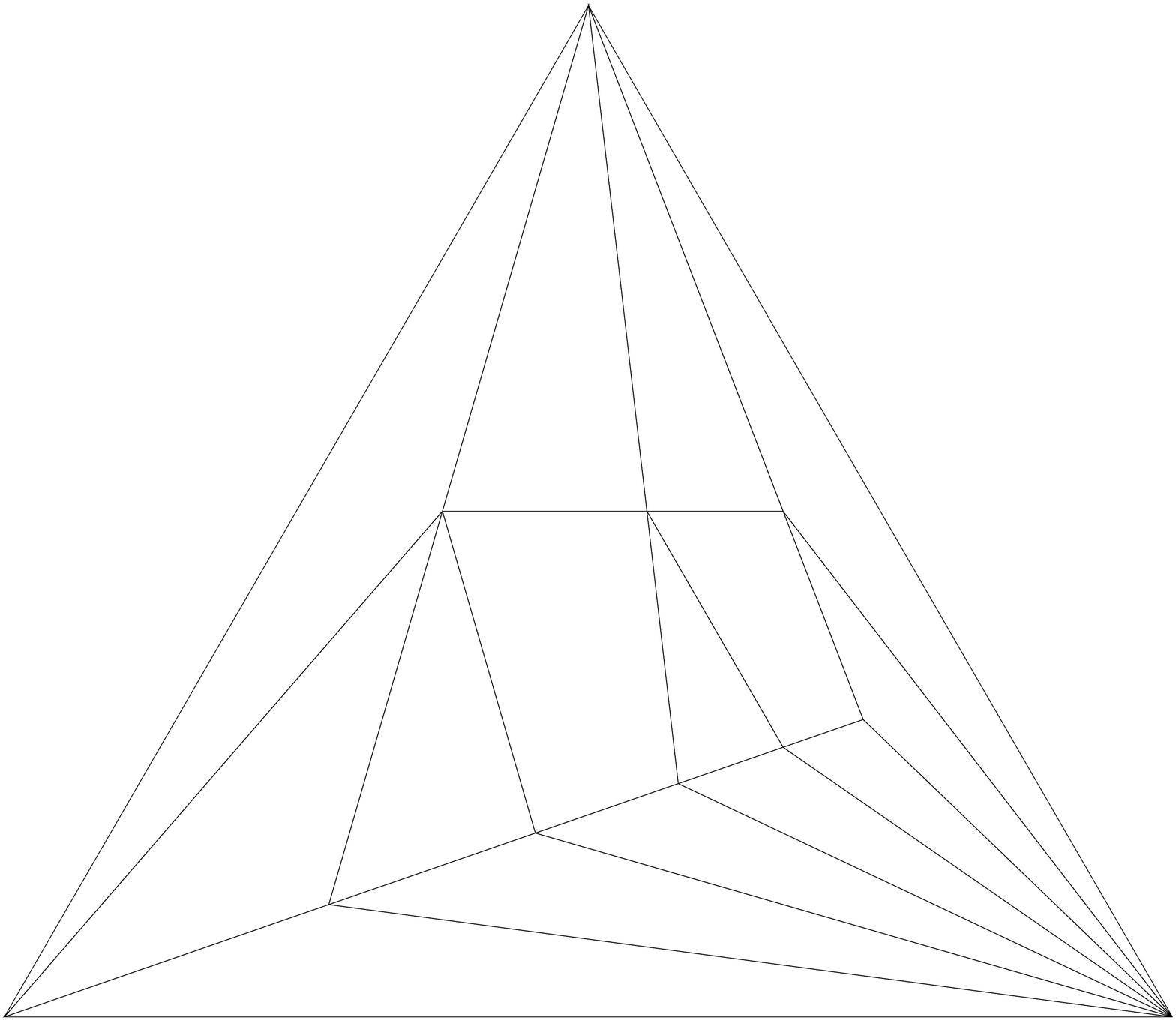,height=3.4cm} 
\end{center}
\caption{ The Gr\"obner fans in the proof of 
Proposition~\ref{prop:notrefinement} intersected with the simplex with 
coordinates $(0,1,0,0)$ (right), $(0,0,1,0)$ (left) and $(0,0,0,1)$ (top).  
The circuit fan is to the left and the toric fan is in the middle.  
The hypergeometric fan is shown to the right.}
\label{fig:toriccircuitfan}
\end{figure}

\begin{proposition}
\label{prop:notrefinement}
In general, neither is the Gr\"obner fan of $\ia$ a refinement of the
Gr\"obner fan of $\ica$, nor vice-versa.
\end{proposition}
\begin{proof}
  Let $A=(7 \; 9 \; 13 \; 15)$.  It is easy to check that
  $\tin_{(0,16,27,1)}(\ica)=\tin_{(0,20,25,3)}(\ica)$ while
  $\tin_{(0,16,27,1)}(\ia)\not=\tin_{(0,20,25,3)}(\ia)$.  Hence
  $(0,16,27,1)$ and $(0,20,25,3)$ lie in the same maximal cell of the
  Gr\"obner fan of $\ica$ but in different maximal cells of the
  Gr\"obner fan of $\ia$. This proves that the Gr\"obner fan of $\ica$
  does not refine the Gr\"obner fan of $\ia$.  On the other hand,
  $\tin_{(0,4,19,9)}(\ia)=\tin_{(0,4,16,5)}(\ia)$ and
  $\tin_{(0,4,19,9)}(\ica)\not=\tin_{(0,4,16,5)}(\ica)$.  Hence the
  Gr\"obner fan of $\ia$ does not refine the Gr\"obner fan of $\ica$.
\end{proof}

The example in Proposition~\ref{prop:notrefinement} would be best
illustrated by a picture of the three-dimensional standard simplex in
$\R^n$ intersected with the fans.  Unfortunately we are limited to two
dimensions in our drawing (Figure \ref{fig:toriccircuitfan}).  The
hypergeometric fan of $\ia$ is drawn at the end of the two Gr\"obner
fans.  

\begin{corollary}  \label{cor:codim2}
If $n-d=2$ the Gr\"obner fan of $\ica$ is a refinement of the Gr\"obner 
fan of $\ia$. 
\end{corollary}

\begin{proof}
  Theorem 3.3.8 in \cite{SST} says that if $n-d=2$ then the Gr\"obner
  fan of $\ia$ equals the hypergeometric fan of $\A$. The corollary
  then follows from Proposition~\ref{prop:top} and the fact that the
  Gr\"obner fan of $\ica$ refines the hypergeometric fan of $\A$.
\end{proof}

\bibliographystyle{plain} 

\bibliography{paper}

\begin{thebibliography}{10}

\bibitem{AtMc}
M.F. Atiyah and I.G. McDonald.
\newblock {\em Introduction to Commutative Algebra}.
\newblock Westview Press, Boulder, CO, 1969.

\bibitem{BFS}
L.~J. Billera, P.~Filliman, and B.~Sturmfels.
\newblock Constructions and complexity of secondary polytopes.
\newblock {\em Advances in Mathematics}, 83:155--179, 1990.

\bibitem{CLO}
D.~Cox, J.~Little, and D.~O'Shea.
\newblock {\em Ideals, Varieties and Algorithms}.
\newblock Springer-Verlag, New York, 1992.

\bibitem{DES}
P.~Diaconis, D.~Eisenbud, and B.~Sturmfels.
\newblock Lattice walks and primary decomposition.
\newblock In B.~E. Sagan and R.~P. Stanley, editors, {\em Mathematical Essays
  in Honor of Gian-Carlo Rota}, pages 173--193, Boston, 1998. {Birkh\"auser}.

\bibitem{DS98}
P.~Diaconis and B.~Sturmfels.
\newblock Algebraic algorithms for sampling from conditional distributions.
\newblock {\em Annals of Statistics}, 26:363--397, 1998.

\bibitem{EiSt}
D.~Eisenbud and B.~Sturmfels.
\newblock Binomial ideals.
\newblock {\em Duke Math. J.}, pages 1--45, 1996.

\bibitem{Ewald}
G.~Ewald.
\newblock {\em Combinatorial Convexity and Algebraic Geometry}, volume 168 of
  {\em Graduate Texts in Mathematics}.
\newblock Springer-Verlag, New York, 1996.

\bibitem{Ful}
W.~Fulton.
\newblock {\em Introduction to Toric Varieties}.
\newblock Princeton University Press, 1993.

\bibitem{M2}
D.~Grayson and M.~Stillman.
\newblock Macaulay 2, a software system for research in algebraic geometry.
\newblock Available at http://www.math.uiuc.edu/Macaulay2.

\bibitem{HS}
S.~Ho\c{s}ten and B.~Sturmfels.
\newblock {G}{R}{I}{N} : An implementation of {G}r{\"o}bner bases for integer
  programming.
\newblock In E.~Balas and J.~Clausen, editors, {\em Integer Programming and
  Combinatorial Optimization}, pages 267--276. Lecture Notes in Computer
  Science 920, Springer Verlag, 1995.

\bibitem{HT2}
S.~Ho\c{s}ten and R.R. Thomas.
\newblock The associated primes of initial ideals of lattice ideals.
\newblock {\em Mathematical Research Letters}, 6:83--97, 1999.

\bibitem{HT01}
S.~Ho\c{s}ten and R.R. Thomas.
\newblock Gomory integer programs.
\newblock {\em Math. Programming, Series B}, 96:271--292, 2003.

\bibitem{Gfan}
A.~Jensen.
\newblock Gfan, a software system for {Gr\"obner} fans.
\newblock Available at http://home.imf.au.dk/ajensen/software/gfan/gfan.html.

\bibitem{MR}
T.~Mora and L.~Robbiano.
\newblock The {G}r{\"o}bner fan of an ideal.
\newblock {\em Journal of Symbolic Computation}, 6:183--208, 1988.

\bibitem{Ohs}
H.~Ohsugi.
\newblock Toric ideals and an infinite family of normal (0,1)-polytopes without
  unimodular regular triangulations.
\newblock {\em Discrete and Computational Geometry}, 27:551--565, 2002.

\bibitem{SST}
M.~Saito, B.~Sturmfels, and N.~Takayama.
\newblock {\em Gr{\"o}bner Deformations of Hypergeometric Differential
  Equations}, volume~6.
\newblock Algorithms and Computation in Mathematics, Springer-Verlag, New-York,
  1999.

\bibitem{Sch}
A.~Schrijver.
\newblock {\em Theory of Linear and Integer Programming}.
\newblock Wiley-Interscience Series in Discrete Mathematics and Optimization,
  New York, 1986.

\bibitem{GBCP}
B.~Sturmfels.
\newblock {\em Gr\"obner Bases and Convex Polytopes}, volume~8 of {\em
  University Lecture Series}.
\newblock American Mathematical Society, Providence, RI, 1996.

\bibitem{STV}
B.~Sturmfels, N.~Trung, and W.~Vogel.
\newblock Bounds on projective schemes.
\newblock {\em Mathematische Annalen}, 302:417--432, 1995.

\end{thebibliography}

\end{document}